\documentclass[11pt]{article}

\usepackage{latexsym}
\usepackage{amsfonts}
\usepackage{amssymb}
\usepackage{mathrsfs}
\usepackage[all]{xy}
\usepackage{epic}
\usepackage{eepic}

\newcommand{\bldd}{\mathbf{D}}
\newcommand{\bldr}{\mathbf{r}}
\newcommand{\blds}{\mathbf{s}}
\newcommand{\bldt}{\mathbf{t}}
\newcommand{\bldm}{\mathbf{m}}
\newcommand{\bldn}{\mathbf{n}}

\newcommand{\sgn}{\mbox{sign\,}}

\newcommand{\be}{\begin{equation}}
\newcommand{\ee}{\end{equation}}
\newcommand{\bea}{\begin{eqnarray}}
\newcommand{\eea}{\end{eqnarray}}
\newcommand{\bean}{\begin{eqnarray*}}
\newcommand{\eean}{\end{eqnarray*}}
\newcommand{\brray}{\begin{array}}
\newcommand{\erray}{\end{array}}

\newcommand{\newsection}[1]{\setcounter{equation}{0}
\setcounter{dfn}{0}
\section{#1}}

\newtheorem{dfn}{Definition}[section]
\newtheorem{thm}[dfn]{Theorem}
\newtheorem{lmma}[dfn]{Lemma}
\newtheorem{ppsn}[dfn]{Proposition}
\newtheorem{crlre}[dfn]{Corollary}
\newtheorem{xmpl}[dfn]{Example}
\newtheorem{rmrk}[dfn]{Remark}

\newcommand{\bdfn}{\begin{dfn}\rm}
\newcommand{\bthm}{\begin{thm}}
\newcommand{\blmma}{\begin{lmma}}
\newcommand{\bppsn}{\begin{ppsn}}
\newcommand{\bcrlre}{\begin{crlre}}
\newcommand{\bxmpl}{\begin{xmpl}}
\newcommand{\brmrk}{\begin{rmrk}\rm}

\newcommand{\edfn}{\end{dfn}}
\newcommand{\ethm}{\end{thm}}
\newcommand{\elmma}{\end{lmma}}
\newcommand{\eppsn}{\end{ppsn}}
\newcommand{\ecrlre}{\end{crlre}}
\newcommand{\exmpl}{\end{xmpl}}
\newcommand{\ermrk}{\end{rmrk}}

\newcommand{\bbc}{\mathbb{C}}

\newcommand{\bbn}{\mathbb{N}}

\newcommand{\scrc}{\mathscr{C}}
\newcommand{\scrf}{\mathscr{F}}

\newcommand{\cla}{\mathcal{A}}

\newcommand{\clh}{\mathcal{H}}

\newcommand{\clk}{\mathcal{K}}
\newcommand{\cll}{\mathcal{L}}

\newcommand{\clg}{\mathcal{G}}

\newcommand{\prf}{\noindent{\it Proof\/}: }

\newcommand{\seq}{\subseteq}

\newcommand{\one}{{1\!\!1}}

\newcommand{\id}{\mbox{id}}

\def \qed { \mbox{}\hfill
$\Box$\vspace{1ex}}

\newcommand{\half}{\frac{1}{2}}

\begin{document}

%%%%%%%%%%%%%%%%%%%%%%%%%%%%%%%%%

%%%%%%%%%%%%%%%%%%%%%%%%%%%%%%%%%
\author{{\sc Partha Sarathi Chakraborty} and
{\sc Arupkumar Pal}}
\title{Equivariant spectral triples
for $SU_q(\ell+1)$ and the odd dimensional
quantum spheres}
\maketitle
%%%%%%%%%%%%%%%%%%%%%%%%%%%%%%%%%%
%%%%%  ABSTRACT
%%%%%%%%%%%%%%%%%%%%%%%%%%%%%%%%%%
 \begin{abstract}
 We formulate the notion of equivariance of an operator
 with respect to a covariant representation of a $C^*$-dynamical system.
    We then use a combinatorial technique used by the authors earlier
   in characterizing spectral triples for $SU_q(2)$ to investigate
   equivariant spectral triples for two classes of spaces:
   the quantum groups $SU_q(\ell+1)$ for $\ell>1$, and
   the odd dimensional quantum spheres $S_q^{2\ell+1}$
   of Vaksman \& Soibelman. In the former case,
   a precise characterization of the sign and the singular
   values of an
   equivariant Dirac operator acting on the $L_2$ space is
   obtained. Using this, we then exhibit
   equivariant Dirac operators with nontrivial sign
   on direct sums of multiple copies of the $L_2$ space.
   In the latter case, viewing $S_q^{2\ell+1}$ as a homogeneous space
   for $SU_q(\ell+1)$, we give a complete characterization
   of equivariant Dirac operators, and also produce an
   optimal family of spectral triples with nontrivial
   $K$-homology class.
 \end{abstract}
{\bf AMS Subject Classification No.:} {\large 58}B{\large 34}, {\large
46}L{\large 87}, {\large
  19}K{\large 33}

%%%%%%%%%%%%%%%%%%%%%%%%%%%%%%%%%%%%%%%%%%%%%%%%%%%

% \tableofcontents

%%%%%%%%%%%%%%%%%%%%%%%%%%%%%%%%%%%%%%%%%%%%%%%%%%%
\newsection{Introduction}
Groups have always played
a very crucial role in the study of geometry of a space, mainly as
objects that govern the symmetry of the space. One would
expect the same in noncommutative geometry also.
Moreover, since one now deals with a larger class of
spaces, mainly noncommutative ones, it is natural to expect that one
would require a larger class, Hopf algebras or the quantum groups, to
play a similar role.  In the classical case, groups which govern
symmetry are themselves nice geometric objects. Here we want to look
at quantum groups from the same angle.
In a previous paper~(\cite{c-p1}), the authors treated the case of the quantum $SU(2)$
group and found a   family of
spectral triples acting on its $L_2$-space
that are equivariant with respect to its natural (co)action.
This family is optimal, in the sense
that given any nontrivial equivariant Dirac operator $D$ acting
on the $L_2$ space, there exists a Dirac  operator $\widetilde{D}$
belonging to this family such that
$\mbox{sign\,}D$ is a compact perturbation of $\mbox{sign\,}\widetilde{D}$
and there exist reals $a$ and  $b$ such that
\[
|D| \leq a + b|\widetilde{D}|.
\]
A generic triple from this family,
that is also a generator of the $K$-homology group,
was analysed by Connes in \cite{co3}
where he used the general theory developed by him and Moscovici
(\cite{c-m}) to  make elaborate
computations and finally ended up with a local index formula.
One beautiful and somewhat surprising observation
in his paper was that the description of the cocycle
given by the difference between the character of the
triple and the cocycle for which index formula was given
involved the Dedekind eta function.
This gave further impetus to the construction of
 spectral triples for quantum groups and their homogeneous
  spaces (\cite{d-s}, \cite{d-l-p-s}, \cite{d-l-s-s-v},
   \cite{h-l},   \cite{kr}, \cite{s-d-l-s-v}).
  It should perhaps be pointed out here
  that the construction by
  Kr\"{a}hmer~(\cite{kr}) is  algebraic in nature and does not
  address the crucial analytic issues involved in the
   definition of a spectral triple.
   The construction by Hawkins \& Landi~(\cite{h-l}) on the other
   hand does not deal with equivariance; and more crucially,
   they restrict themselves to the construction of bounded
   Kasparov modules. But in Noncommutative geometry,
   spectral triples or the unbounded Kasparov modules
   are key ingredients, as they work as a looking glass allowing
   one to distinguish between continuous and smooth functions.

Our aim in the present paper is to look for higher
dimensional counterparts of the spectral triples
found in \cite{c-p1}.
We first formulate precisely what
one means by an equivariant spectral triple in a general
set up (this is already implicit in \cite{c-p1})
and then study equivariant Dirac operators
for two classes of spaces, both of which can
be thaught of as higher dimensional analogues
of  $SU_q(2)$ which was worked out earlier.
First, we analyse equivariant Dirac operators
acting on the $L_2$-spaces of the groups $SU_q(\ell+1)$.
We derive a precise expression for the singular
values of an equivariant Dirac operator, and show that
a Dirac operator with these singular values will
have the correct summability property.
We also show that for $\ell>1$,
an equivariant Dirac operator  acting on $L_2(G)$
have to have trivial sign.
Thus for $\ell>1$, one would be forced to bring in multiplicity
when looking for equivariant Dirac operators
with nontrivial sign.
Using this observation, we then exhibit a
family of equivariant Dirac operators
acting on direct sums of multiple copies of
the $L_2$ space and having nontrivial sign.
Whether these Dirac operators
have nontrivial $K$-homology class is still not known.
In the last section, we take up the odd dimensional
quantum spheres $S_q^{2\ell+1}$. In this case, the outcome
turns out to be more satisfactory. After characterizing
the sign and the singular values of Dirac operators on
$L_2(S_q^{2\ell+1})$ equivariant under the action of
the group $SU_q(\ell+1)$, we produce, just like
in the $SU_q(2)$ case, an optimum family of
nontrivial equivariant Dirac operators that are
$(2\ell+1)$-summable.

The paper is organised as follows.
In the next section, we will recall from~\cite{c-p0} the
combinatorial method that was earlier used implicitly
in \cite{c-p1} and \cite{c-p2}.
In section~3, we formulate the notion of equivariance.
This has been done using the quantum group at
the function algebra level rather than passing on
to the quantum universal envelopping algebra level.
In section~4, we briefly recall the quantum group
$SU_q(\ell+1)$ and its representation theory.
In particular, we describe a nice basis for
the $L_2$ space and study the Clebsch-Gordon
coefficients. These are used in section~5 to describe the
action by left multiplication on the $L_2$ space explicitly.
In section~6, we write down the conditions
coming from the boundedness of commutators with $D$.
In sections~7 and 8, we analyze the equivariant Dirac operators
for $SU_q(\ell+1)$. First we give a precise characterization of
the singular values in section~7,
and then a characterization of the sign in section~8.
In section~9, we deal with the odd dimensional quantum spheres.

%%%%%%%%%%%%%%%%%%%%%%%%%%%%%%%%%%%%%%%%%%%%%%%%%%%
\newsection{The general scheme}
%%%%%%%%%%%%%%%%%%%%%%%%%%%%%%%%%%%%%%%%%%%%%%%%%%%
Let us recall the combinatorial set up from~\cite{c-p0}.

Suppose $\clh$ is a Hilbert space,
and $D$ is a self-adjoint operator on $\clh$ with compact resolvent.
Then $D$ admits a spectral resolution $\sum_{\gamma\in\Gamma} d_\gamma
P_\gamma$, where the $d_\gamma$'s are all distinct and each $P_\gamma$
is a finite dimensional projection.  Assume now onward that all the
$d_\gamma$'s are nonzero. Let $c$ be a positive real.  Let
us  define a graph $\clg_c$ as follows: take the vertex set $V$ to
be $\Gamma$.  Connect two vertices $\gamma$ and $\gamma'$ by an edge
if $|d_\gamma-d_{\gamma'}|<c$. Let $V^+=\{\gamma\in V: d_\gamma>0\}$ and
$V^-=\{\gamma\in V: d_\gamma<0\}$.  This will give us a partition of
$V$.
This partition has the following important property:
there does not exist infinite number of disjoint
paths each going from a point in $V^+$
to a point in $V^-$. Here disjoint paths mean paths
for which the set of vertices of one does not
intersect the set of vertices
of the other.
This is easy to see, because
if there is a path from $\gamma$ to $\delta$
and $d_\gamma>0$, $d_\delta<0$, then for some $\alpha$ on the path,
one must have $d_\alpha\in[-c,c]$.
Since the paths are disjoint, it
would contradict the compact resolvent condition.
We will call such a partition a sign-determining partition.

We will use this knowledge about the graph.
We start with an equivariant operator that is self-adjoint
and has discrete spectrum. Equivariance will give us
an idea about the spectral resolution
$\sum_{\gamma\in\Gamma}d_\gamma P_\gamma$.
Next we use the action of the algebra elements on the basis
elements of $\clh$ and the boundedness of their
commutators with $D$.
This gives certain growth restrictions
on the $d_\gamma$'s. These will give us some information
about the edges in the graph. We exploit this knowledge
to characterize those partitions $(V_1,V_2)$ of the vertex set
that are sign-determining, i.\ e.\ do not admit any infinite ladder.
The sign of the operator $D$ must be of the form
$\sum_{\gamma\in V_1}P_\gamma-\sum_{\gamma\in V_2}P_\gamma$
where  $(V_1,V_2)$ is a sign-determining partition.
  Of course, for a given $c$, the graph
$\clg_c$ may have no edges, or too few edges (if the singular values
of $D$ happen to grow too fast), in which case, we will
be left with too many sign-determining
partitions.
Fortunately, the operators we are interested in are meant to be the
Dirac operators of some commutative/noncommutative manifold. Therefore
the singular values of $D$ will grow at the rate of $O(n^{1/d})$ for
some $d\geq 1$.  So one can choose a large enough $c$ and work with
the graph $\clg_c$.
In other words, we would like to characterize
those partitions that are sign-determining
for all sufficiently large values of $c$.

%%%%%%%%%%%%%%%%%%%%%%%%%%%%%%%%%%%%%%%%%%%%%%%%%%%%%%%%%
%%%%%%%%%%%%%%%%%%%%%%%%%%%%%%%%%%%%%%%%%%%%%%%%%%%%%%%%%
%%%%%%%%%%%%%%%%%%%%%%%%%%%%%%%%%%%%%%%%%%%%%%%%%%%%%%%%%
%%%%%%%%%%%%%%%%%%%%%%%%%%%%%%%%%%%%%%%%%%%%%%%%%%%%%%%%%
%%%%%%%%%%%%%%%%%%%%%%%%%%%%%%%%%%%%%%%%%%%%%%%%%%%%%%%%%
%%%%%%%%%%%%%%%%%%%%%%%%%%%%%%%%%%%%%%%%%%%%%%%%%%%%%%%%%
%%%%%%%%%%%%%%%%%%%%%%%%%%%%%%%%%%%%%%%%%%%%%%%%%%%%%%%%%
%%%%%%%%%%%%%%%%%%%%%%%%%%%%%%%%%%%%%%%%%%%%%%%%%%%%%%%%%
%%%%%%%%%%%%%%%%%%%%%%%%%%%%%%%%%%%%%%%%%%%%%%%%%%%%%%%%%
%%%%%%%%%%%%%%%%%%%%%%%%%%%%%%%%%%%%%%%%%%%%%%%%%%%%%%%%%
%%%%%%%%%%%%%%%%%%%%%%%%%%%%%%%%%%%%%%%%%%%%%%%%%%%%%%%%%
%%%%%%%%%%%%%%%%%%%%%%%%%%%%%%%%%%%%%%%%%%%%%%%%%%%%%%%%%

%%%%%%%%%%%%%%%%%%%%%%%%%%%%%%%%%%%%%%%%%%%%%%%%%%%%%%%%%
\newsection{Equivariance}
%%%%%%%%%%%%%%%%%%%%%%%%%%%%%%%%%%%%%%%%%%%%%%%%%%%
Suppose $G$ is a compact group, quantum or classical,
and $\cla$ is a unital $C^*$-algebra. Assume that
$G$ has an action on $\cla$ given by
$\tau:\cla\rightarrow\cla\otimes C(G)$,
so that $(\id\otimes\Delta)\tau=(\tau\otimes\id)\tau$,
$\Delta$ being the coproduct.
In other words, we have a $C^*$-dynamical system $(\cla,G,\tau)$.
Our goal is to study spectral triples for $\cla$
equivariant under this action.
Let us first say what we mean by `equivariant' here.

A covariant representation $(\pi,u)$
of $(\cla,G,\tau)$ consists of
a unital *-representation $\pi:\cla\rightarrow\cll(\clh)$,
a unitary representation $u$ of $G$ on $\clh$, i.e.\
a unitary element of the multiplier algebra $M(\clk(\clh)\otimes C(G))$
such that they obey the condition
$(\pi\otimes\id)\tau(a)=u(\pi(a)\otimes I)u^*$ for all $a\in\cla$.
%%%%%%%%%%%%%%%%%%%%%%%%%%%%%%%%%%%%%%%%%%%%%
\bdfn
Suppose $(\cla, G,\tau)$ is a $C^*$-dynamical system.
An operator $D$ acting on a Hilbert space $\clh$
is said to be \textbf{equivariant} with respect to a covariant
representation $(\pi,u)$ of the system if
$D\otimes I$ commutes with $u$.
\edfn
%%%%%%%%%%%%%%%%%%%%%%%%%%%%%%%%%%%%%%%%%%%%%

Since the operator $D$ is self-adjoint with compact resolvent, it will
admit a spectral resolution $\sum_\lambda d_\lambda P_\lambda$, where
the $d_\lambda$'s are distinct and each $P_\lambda$ is finite
dimensional.  Also, $D$ has been assumed to be equivariant --- so that
the $P_\lambda$'s commute with $u$ (to be precise, the
$(P_\lambda\otimes I)$'s do), i.e.\ $u$ keeps each $P_\lambda\clh$
invariant. As $G$ is compact, each $P_\lambda\clh$ will decompose
further as $\oplus_\mu P_{\lambda\mu}\clh$ such that the restriction
of $u$ to each $P_{\lambda\mu}$ is irreducible.  In other words, one
can now write $D$ in the form $\sum_{\gamma\in\Gamma}d_\gamma
P_\gamma$ for some index set $\Gamma$ and a family of finite
dimensional projections $P_\gamma$ such that each $P_\gamma$ commutes
with $u$ and the restriction of $u$ to each $P_\gamma$ is irreducible.

In this paper, we will deal with two cases,
the group in question in both cases will be $G=SU_q(\ell+1)$.
The $C^*$-algebra $\cla$ on which the group
acts will be $C(SU_q(\ell+1))$ in one case and
$C(S_q^{2\ell+1})$ in the other.
Let us discuss the first case a little here.
The action $\tau$ here will be the natural action coming from the
coproduct, $\clh$ is $L_2(G)$, $\pi$ is the
representation of $\cla=C(SU_q(\ell+1))$
on $\clh$ by left multiplication, and $u$ is the right regular
representation.  Structure of the regular representation of a compact
(quantum) group along with the remarks made above tell us the
following.  Let $\Lambda$ be the set of unitary irreducible
representation-types for $G$. Then $\clh$ decomposes as
$\oplus_{\lambda\in\Lambda}\clh_\lambda$, where the restriction of $u$
to $\clh_\lambda$ is equivalent to $\mbox{dim}\,\lambda$ copies of the
irreducible $\lambda$, and also that $D$ respects this decomposition.
Further, restriction of $D$ to $\clh_\lambda$ is of the form
$\sum_{\mu}d_{\lambda\mu}P_{\lambda\mu}$, $u$ commutes with each of
these $P_{\lambda\mu}$'s, and the restriction of $u$ to
$P_{\lambda\mu}\clh$ is equivalent to $\lambda$. Let $N_\lambda$ be
any set with $|N_\lambda|=\mbox{dim}\,\lambda$.  One can then choose
an orthonormal basis $\{e^\lambda_{ij}:i,j\in N_\lambda\}$ such that
the spaces $P_{\lambda\mu}\clh$ are precisely
$\mbox{span}\,\{e^\lambda_{ij}:j\in N_\lambda\}$ for distinct values
of $i\in N_\lambda$.  Since $D$ is of the form $\sum_\lambda\sum_\mu
d_{\lambda\mu}P_{\lambda\mu}$, in this system of bases, $D$ will look
like $e^\lambda_{ij}\mapsto d(\lambda,i)e^\lambda_{ij}$.  In what
follows, we will make a special choice of $N_\lambda$, which will make
the combinatorial analysis very convenient.

%%%%%%%%%%%%%%%%%%%%%%%%%%%%%%%%%%%%%%%%%%%%%%%%%
\newsection{Preliminaries on $SU_q(\ell+1)$}
%%%%%%%%%%%%%%%%%%%%%%%%%%%%%%%%%%%%%%%%%%%%%%%%%
Let $\mathfrak{g}$ be a complex simple Lie algebra of rank $\ell$.
let $(\!(a_{ij})\!)$ be the associated Cartan matrix,
$q$ be a real number lying in the interval $(0,1)$
and let $q_i=q^{(\alpha_i,\alpha_i)/2}$, where $\alpha_i$'s are the simple roots
of $\mathfrak{g}$.
Then the quantised universal envelopping algebra (QUEA)
$U_q(\mathfrak{g})$ is the algebra
generated by $E_i$, $F_i$, $K_i$ and $K_i^{-1}$, $i=1,\ldots,\ell$, satisfying the
following relations
\begin{displaymath}
K_iK_j=K_jK_i,\quad K_iK_i^{-1}=K_i^{-1}K_i=1,
\end{displaymath}
\begin{displaymath}
K_iE_jK_i^{-1}=q_i^{\half a_{ij}}E_j,\quad
K_iF_jK_i^{-1}=q_i^{-\half a_{ij}}F_j,
\end{displaymath}
\begin{displaymath}
E_iF_j-F_jE_i=\delta_{ij}\frac{K_i^2-K_i^{-2}}{q_i-q_i^{-1}},
\end{displaymath}
\begin{displaymath}
\sum_{r=0}^{1-a_{ij}}(-1)^r{{1-a_{ij}}\choose r}_{q_i}
  E_i^{1-a_{ij}-r}E_jE_i^r =0 \quad\forall\, i\neq j,
\end{displaymath}
\begin{displaymath}
\sum_{r=0}^{1-a_{ij}}(-1)^r{{1-a_{ij}}\choose r}_{q_i}
  F_i^{1-a_{ij}-r}F_jF_i^r =0\quad \forall\, i\neq j,
\end{displaymath}
where ${n\choose r}_q$ denote the $q$-binomial coefficients.
Hopf *-structure comes from the following maps:
\[
\Delta(K_i)=K_i\otimes K_i,\quad \Delta(K_i^{-1})=K_i^{-1}\otimes K_i^{-1},
\]
\[
\Delta(E_i)=E_i\otimes K_i  + K_i^{-1}\otimes E_i,\quad
\Delta(F_i)=F_i\otimes K_i  + K_i^{-1}\otimes F_i,
\]
\[
\epsilon(K_i)=1,\quad \epsilon(E_i)=0=\epsilon(F_i),
\]
\[
S((K_i)=K_i^{-1},\quad S(E_i)=-q_iE_i,\quad S(F_i)=-q_i^{-1}F_i,
\]
\[
K_i^*=K_i,\quad E_i^*=-q_i^{-1}F_i,\quad F_i^*=-q_iE_i.
\]

In the type A case, the associated Cartan matrix is given by
\[
a_{ij}=\cases{2& if $i=j$,\cr
              -1 & if $i=j\pm1$,\cr
               0 & otherwise,}
\]
and $(\alpha_i,\alpha_i)=2$ so that $q_i=q$ for all $i$.
The QUEA in this case is denoted by $u_q(su(\ell+1))$.

Take the collection of matrix entries of all finite-dimensional
unitarizable $u_q(su(\ell+1))$-modules. The algebra generated by these
gets a natural Hopf*-structure as the dual of $u_q(su(\ell+1))$.  One
can also put a natural $C^*$-norm on this.  Upon completion with
respect to this norm, one gets a unital $C^*$-algebra that plays the
role of the algebra of continuous functions on $SU_q(\ell+1)$.  For a
detailed account of this, refer to chapter~3, \cite{ko-so}.  In
\cite{w}, Woronowicz gave a different description of this
$C^*$-algebra.  which was later shown by Rosso (\cite{r}) to be
equivalent to the earlier one.

For remainder of this article, we will take $G$ to be $SU_q(\ell+1)$
and $\cla$ will be the $C^*$-algebra of continuous functions on $G$.

%%%%%%%%%%%%%%%%%%%%%%%%%%%%%%%%%%%%%%%%%%%%%%%%%%%
\paragraph{Gelfand-Tsetlin tableaux.}
%%%%%%%%%%%%%%%%%%%%%%%%%%%%%%%%%%%%%%%%%%%%%%%%%%%
Irreducible unitary representations of the group
$SU_q(\ell+1)$ are indexed by
Young tableaux $\lambda=(\lambda_1,\ldots,\lambda_{\ell+1})$,
where $\lambda_i$'s are nonnegative integers,
$\lambda_1\geq \lambda_2\geq \ldots\geq \lambda_{\ell+1}$
(Theorem~1.5, \cite{w}).
Write $\clh_\lambda$ for the Hilbert space where
the irreducible $\lambda$ acts.
There are various ways of indexing the basis elements
of $\clh_\lambda$. The one we will use is due to Gelfand
and Tsetlin.
According to their prescription, basis elements for
$\clh_\lambda$ are parametrized by arrays of the form
\[
\bldr=\left(\matrix{r_{11}&r_{12} &\cdots&r_{1,\ell}&r_{1,\ell+1}\cr
                      r_{21}&r_{22}&\cdots &r_{2,\ell}&\cr
                          &\cdots&&&\cr
                      r_{\ell,1}&r_{\ell,2}&&&\cr
                      r_{\ell+1,1}&&&&}\right),
\]
where $r_{ij}$'s are integers satisfying
$r_{1j}=\lambda_j$ for $j=1,\ldots,\ell+1$,
$r_{ij}\geq r_{i+1,j}\geq r_{i,j+1}\geq 0$ for all $i$, $j$.
Such arrays are known as Gelfand-Tsetlin tableaux, to be abreviated
as GT tableaux for the rest of this section.
For a GT tableaux $\bldr$, the symbol $\bldr_{i\cdot}$ will denote its
$i$\raisebox{.4ex}{th} row.
It is well-known that two representations indexed respectively
by $\lambda$ and $\lambda'$ are equivalent if and only if
$\lambda_j-\lambda_j^\prime$ is independent of $j$ (\cite{w}).
Thus one gets an equivalence relation on the set of Young tableaux
$\{ \lambda=(\lambda_1,\ldots,\lambda_{\ell+1}):
\lambda_1\geq \lambda_2\geq \ldots\geq \lambda_{\ell+1}, \lambda_j\in\bbn\}$.
This, in turn, induces an equivalence relation on the set of
all GT tableaux $\Gamma=\{\bldr: r_{ij}\in\bbn,
  r_{ij}\geq r_{i+1,j}\geq r_{i,j+1}\}$: one says $\bldr$ and $\blds$
are equivalent if $r_{ij}-s_{ij}$ is independent of $i$ and $j$.
By $\Gamma$ we will mean the above set modulo this equivalence.

We will denote by $u^\lambda$ the irreducible unitary indexed by $\lambda$,
$\{e(\lambda,\bldr):\bldr_{1\cdot}=\lambda\}$ will denote an orthonormal basis
for $\clh_\lambda$ and $u^\lambda_{\bldr\blds}$ will stand for the matrix entries
of $u^\lambda$ in this basis. The symbol $\one$ will denote the Young tableaux
$(1,0,\ldots,0)$. We will often omit the symbol $\one$
and just write $u$ in order to denote $u^\one$.
Notice that any GT tableaux $\bldr$ with first row $\one$
must be, for some $i\in\{1,2,\ldots,\ell+1\}$, of the form $(r_{ab})$, where
\[
r_{ab}=\cases{1 &if $1\leq a\leq i$ and $b=1$,\cr
       0   &otherwise.}
\]
Thus such a GT tableaux is uniquely determined by the integer $i$.
We will write just $i$ for this GT tableaux $\bldr$.
Thus for example, a typical matrix entry of $u^\one$ will be
written simply as $u_{ij}$.

Let $\bldr=(r_{ab})$ be a GT tableaux.
Let
$H_{ab}(\bldr):=r_{a+1,b}-r_{a,b+1}$ and
$V_{ab}(\bldr):=r_{ab}-r_{a+1,b}$.
An element $\bldr$ of $\Gamma$ is completely
specified by the following differences
\[
\bldd(\bldr)=\left(\matrix{V_{11}(\bldr)&H_{11}(\bldr)
         &H_{12}(\bldr)&\cdots&H_{1,\ell-1}(\bldr)&H_{1,\ell}(\bldr)\cr
  V_{21}(\bldr)&H_{21}(\bldr)&H_{22}(\bldr)&\cdots&H_{2,\ell-1}(\bldr)&\cr
         &\cdots&&&&\cr
       V_{\ell,1}(\bldr)&H_{\ell,1}(\bldr)&&&&}\right).
\]
The differences satisfy the following inequalities
\be\label{ineq}
\sum_{k=0}^b H_{a-k,k+1}(\bldr)\leq V_{a+1,1}(\bldr)
       +\sum_{k=0}^b H_{a-k+1,k+1}(\bldr),\quad
1\leq a\leq \ell,\;\;0\leq b\leq a-1.
\ee
Conversely, if one has an array of the form
\[
\left(\matrix{V_{11}&H_{11}&H_{12}&\cdots&H_{1,\ell-1}&H_{1,\ell}\cr
       V_{21}&H_{21}&H_{22}&\cdots&H_{2,\ell-1}&\cr
         &\cdots&&&&\cr
       V_{\ell,1}&H_{\ell,1}&&&&}\right),
\]
where $V_{ij}$'s and $H_{ij}$'s are in $\bbn$ and obey
the inequalities~(\ref{ineq}), then the above array is of the form
$\bldd(\bldr)$ for some GT tableaux $\bldr$. Thus the quantities
$V_{a1}$ and $H_{ab}$ give a coordinate system for elements in $\Gamma$.
The following diagram explains this new coordinate system.
The hollow circles stand for the $r_{ij}$'s.
The entries are decreasing along the direction of the arrows,
and the $V_{ij}$'s and the $H_{ij}$'s are the difference
between the two endpoints of the corresponding arrows.\\

%%%%%%%%%%%%%%%%%%%%%%%%%%%%%%%%%%%%%%%%%%%%%%%%
\hspace*{100pt}
  \def\labelstyle{\scriptstyle}
 
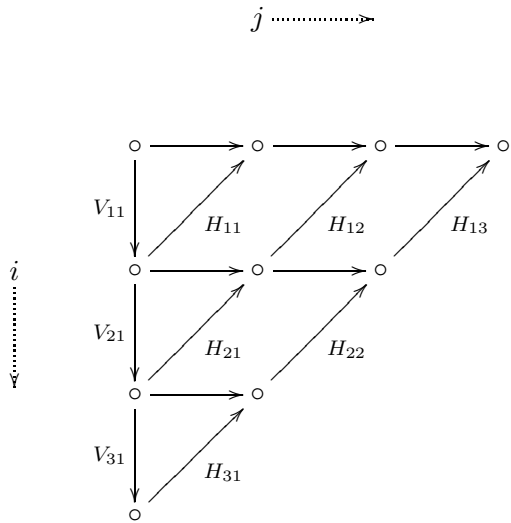
\xymatrix@C=35pt@R=35pt{
   &  &  j\ar@{.>}[r] &&\\
   & \circ\ar@{->}[r]\ar@{->}[d]_{V_{11}} &  \circ\ar@{->}[r]
                       & \circ\ar@{->}[r] &\circ\\
i\ar@{.>}[d] &  \circ\ar@{->}[r]\ar@{->}[d]_{V_{21}}\ar@{->}[ur]_{H_{11}} &
    \circ\ar@{->}[r]\ar@{->}[ur]_{H_{12}} & \circ\ar@{->}[ur]_{H_{13}} & \\
   & \circ\ar@{->}[r]\ar@{->}[d]_{V_{31}}\ar@{->}[ur]_{H_{21}} &
               \circ\ar@{->}[ur]_{H_{22}} &\\
   &  \circ\ar@{->}[ur]_{H_{31}} & }\\
%%%%%%%%%%%%%%%%%%%%%%%%%%%%%%%%%%%%%%%%%%%%%%%%

%%%%%%%%%%%%%%%%%%%%%%%%%%%%%%%%%%%%%%%%%%%%%%%%%%%%%%%%%%%
%% CG coefficients
%%%%%%%%%%%%%%%%%%%%%%%%%%%%%%%%%%%%%%%%%%%%%%%%%%%%%%%%%%%
\paragraph{Clebsch-Gordon coefficients.}
Look at the representation $u^\one\otimes u^\lambda$
acting on $\clh_\one\otimes\clh_\lambda$.
The representation decomposes as a direct sum
$\oplus_\mu u^\mu$, i.e.\ one has a corresponding
decomposition $\oplus_\mu\clh_\mu$ of $\clh_\one\otimes\clh_\lambda$.
Thus one has two orthonormal bases
$\{e^\mu_\blds\}$ and $\{e^\one_i\otimes e^\lambda_\bldr\}$.
The Clebsch-Gordon coefficient $C_q(\one,\lambda,\mu;i,\bldr,\blds)$
is defined to be the inner product
$\langle e^\mu_\blds, e^\one_i\otimes e^\lambda_\bldr\rangle$.
Since $\one$, $\lambda$ and $\mu$ are just the first rows of
$i$, $\bldr$ and $\blds$ respectively, we will often denote
the above quantity just by $C_q(i,\bldr,\blds)$.

Next, we will compute the quantities $C_q(i,\bldr,\blds)$.  We will
use the calculations given in (\cite{k-s}, pp.\ 220), keeping in mind
that for our case (i.e.\ for $SU_q(\ell+1)$), the top right entry of
the GT tableaux is zero.

Let $M=(m_1,m_2,\ldots,m_i)\in\bbn^i$ be such that $1\leq m_j\leq \ell+2-j$.
Denote by $M(\bldr)$ the tableaux $\blds$ defined by
\be\label{movenotation}
s_{jk}=\cases{r_{jk}+1 & if $k=m_j$, $1\leq j\leq i$,\cr
         r_{jk} & otherwise.}
\ee
With this notation, observe now that
$C_q(i,\bldr,\blds)$ will be zero unless $\blds$ is
$M(\bldr)$ for some $M\in\bbn^i$.
(One has to keep in mind though that not all tableaux of the form $M(\bldr)$
is a valid GT tableaux)

From (\cite{k-s}, pp.\ 220), we have
\be\label{cgc1}
C_q(i,\bldr,M(\bldr))=\prod_{a=1}^{i-1}
\left\langle \brray{ll}
                (1,\mathbf{0}) &\bldr_{a\cdot} \cr
                (1,\mathbf{0}) &\bldr_{a+1\cdot}
             \erray\left| \brray{l}
                   \bldr_{a\cdot}+e_{m_a}\cr
                   \bldr_{a+1\cdot}+e_{m_{a+1}}
              \erray\right.\right\rangle
\times
\left\langle \brray{ll}
                (1,\mathbf{0}) &\bldr_{i\cdot} \cr
                (0,\mathbf{0}) &\bldr_{i+1\cdot}
             \erray\left| \brray{l}
                   \bldr_{i\cdot}+e_{m_i}\cr
                   \bldr_{i+1\cdot}
              \erray\right.\right\rangle,
\ee
where $e_k$ stands for a vector (in the appropriate space) whose
$k$\raisebox{.4ex}{th} coordinate is 1 and the rest are all zero, and
\bea
\left\langle \brray{ll}
                (1,\mathbf{0}) &\bldr_{a\cdot} \cr
                (1,\mathbf{0}) &\bldr_{a+1\cdot}
             \erray\left| \brray{l}
                   \bldr_{a\cdot}+e_j\cr
                   \bldr_{a+1\cdot}+e_k
              \erray\right.\right\rangle^2
&=&
q^{-r_{aj}+r_{a+1,k} - k+j}
\times
\prod_{{i=1}\atop{i\neq j}}^{\ell+2-a}
   \frac{[r_{a,i}-r_{a+1,k}-i+k]_q  }{[r_{a,i}-r_{a,j}-i+j]_q} \nonumber \\
&& \times
\prod_{{i=1}\atop{i\neq k}}^{\ell+1-a}
   \frac{[r_{a+1,i}-r_{a,j}-i+j-1]_q  }{[r_{a+1,i}-r_{a+1,k}-i+k-1]_q},\label{corrected_1}\\
\left\langle \brray{ll}
                (1,\mathbf{0}) &\bldr_{a\cdot} \cr
                (0,\mathbf{0}) &\bldr_{a+1\cdot}
             \erray\left| \brray{l}
                   \bldr_{a\cdot}+e_j\cr
                   \bldr_{a+1\cdot}
              \erray\right.\right\rangle^2
&=& q^{\left(1-j+\sum_{i=1}^{\ell+1-a}r_{a+1,i} -
           \sum_{{i=1}\atop{i\neq j}}^{\ell+2-a}r_{a,i}\right)} \nonumber \\
&&
\times \left(
\frac{\prod_{i=1}^{\ell+1-a}[r_{a+1,i}-r_{aj}-i+j-1]_q  }
{\prod_{{i=1}\atop{i\neq j}}^{\ell+2-a}[r_{a,i}-r_{aj}-i+j]_q  }\right),
\label{corrected_2}
\eea
where for an integer $n$, $[n]_q$ denotes the $q$-number $(q^n-q^{-n})/(q-q^{-1})$.
After some lengthy but straightforward computations,
we get the following two relations:
\be
\left|
\left\langle \brray{ll}
                (1,\mathbf{0}) &\bldr_{a\cdot} \cr
                (1,\mathbf{0}) &\bldr_{a+1\cdot}
             \erray\left| \brray{l}
                   \bldr_{a\cdot}+e_j\cr
                   \bldr_{a+1\cdot}+e_k
              \erray\right.\right\rangle
\right| = A'q^A,
\ee
\be
\left|
\left\langle \brray{ll}
                (1,\mathbf{0}) &\bldr_{a\cdot} \cr
                (0,\mathbf{0}) &\bldr_{a+1\cdot}
             \erray\left| \brray{l}
                   \bldr_{a\cdot}+e_j\cr
                   \bldr_{a+1\cdot}
              \erray\right.\right\rangle
\right| = B'q^B,
\ee
where
\bea
A&=&\cases{\displaystyle{\sum_{j\wedge k < b < j\vee k}(r_{a+1,b}-r_{a,b})}
              +(r_{a+1,j\wedge k}-r_{a,j\vee k})  & if $j\neq k$,\cr
   0 & if $j=k$.} \cr
&=& \sum_{j\wedge k \leq b < j\vee k}(r_{a+1,b}-r_{a,b+1})
    +2 \sum_{k < b < j}(r_{a,b}-r_{a+1,b}) \cr
&=&  \sum_{j\wedge k \leq b < j\vee k}H_{ab}(\bldr)
             + 2 \sum_{k < b < j}V_{ab}(\bldr).\label{cgc2}\\
B &=&  \sum_{j \leq b < \ell+2-a}H_{ab}(\bldr),\label{cgc3}
\eea
and $A'$ and $B'$ both lie between two positive constants
independent of $\bldr$, $a$, $j$ and $k$
(Here and elsewhere in this paper, an empty summation
would always mean zero).

Combining these, one gets
\be \label{cgc4}
C_q(i,\bldr, M(\bldr))=P\cdot q^{C(i,\bldr,M)},
\ee
where
\be \label{cgc5}
C(i,\bldr,M)=\sum_{a=1}^{i-1}\left(
   \sum_{m_a\wedge m_{a+1} \leq b < m_a\vee m_{a+1}}H_{ab}(\bldr)
   +2 \sum_{m_{a+1} < b < m_a}V_{ab}(\bldr)\right)
+\sum_{m_i \leq b < \ell+2-i}H_{ib}(\bldr),
\ee
and $P$ lies between two positive constants
that are independent of $i$, $\bldr$ and $M$.

\brmrk
The formulae (\ref{corrected_1}) and (\ref{corrected_2})
are obtained from equations~(45) and (46), page 220, \cite{k-s}
by replacing $q$ with $q^{-1}$. Equation~(45) is a special
case of the more general formula (48), page 221, \cite{k-s}.
However, there is a small error in equation~(48) there.
The correct form can be found in equations~(3.1, 3.2a, 3.2b)
in \cite{a-s}. That correction has been incorporated in
equations~(\ref{corrected_1}) and (\ref{corrected_2}) here.
\ermrk

%%%%%%%%%%%%%%%%%%%%%%%%%%%%%%%%%%%%%%%%%%%%%%%%%%%%%%%%%%%%%%%%%%%
\newsection{Left multiplication operators}
%%%%%%%%%%%%%%%%%%%%%%%%%%%%%%%%%%%%%%%%%%%%%%%%%%%%%%%%%%%%%%%%%%%
The matrix entries $u^\lambda_{\bldr\blds}$ form a complete orthogonal set
of vectors in $L_2(G)$. Write $e^\lambda_{\bldr\blds}$ for
$\|u^\lambda_{\bldr\blds}\|^{-1}u^\lambda_{\bldr\blds}$.
Then the $e^\lambda_{\bldr\blds}$'s form a complete orthonormal basis
for $L_2(G)$. Let $\pi$ denote the representation of $\cla$ on
$L_2(G)$ by left multiplications. We will now derive an expression for
$\pi(u_{ij})e^\lambda_{\bldr\blds}$.

From the definition of matrix entries and that of the CG coefficients,
one gets
\be \label{cb1}
u^\rho e(\rho,\bldt)=\sum_\blds u^\rho_{\blds\bldt}e(\rho,\blds),
\ee
\be \label{cb2}
e(\mu,\bldn)=\sum_{j,\blds}C_q(j,\blds,\bldn)e(\one,j)\otimes e(\lambda,\blds).
\ee
Apply $u\otimes u^\lambda$ on both sides and note that
$u\otimes u^\lambda$ acts on $e(\mu,\bldn)$ as $u^\mu$:
\be \label{cb3}
\sum_\bldm u^\mu_{\bldm\bldn}e(\mu,\bldm)=
\sum_{j,\blds}\sum_{i,\bldr}C_q(j,\blds,\bldn)
 u_{ij}u^\lambda_{\bldr\blds}e(\one,i)\otimes e(\lambda,\bldr).
\ee
Next, use (\ref{cb2}) to expand $e(\mu,\bldm)$ on the left hand side to get
\be
\sum_{i,\bldr,\bldm} u^\mu_{\bldm\bldn}
   C_q(i,\bldr,\bldm)e(\one,i)\otimes e(\lambda,\bldr)
=
\sum_{j,\blds}\sum_{i,\bldr}C_q(j,\blds,\bldn)
        u_{ij}u^\lambda_{\bldr\blds}e(\one,i)\otimes
         e(\lambda,\bldr).
\ee
Equating coefficients, one gets
\be
\sum_{\bldm} C_q(i,\bldr,\bldm)u^\mu_{\bldm\bldn}
=
\sum_{j,\blds}C_q(j,\blds,\bldn)
        u_{ij}u^\lambda_{\bldr\blds}.
\ee
Now using orthogonality of the matrix
$(\!(C_q(\one,\lambda,\mu;j,\blds,\bldn))\!)_{(\mu,\bldn),(j,\blds)}$,
we obtain
\be\label{alg_left_mult}
u_{ij}u^\lambda_{\bldr\blds}
= \sum_{\mu,\bldm,\bldn}
  C_q(i,\bldr,\bldm)C_q(j,\blds,\bldn)u^\mu_{\bldm\bldn}.
\ee
From (\cite{k-s}, pp.\ 441), one has
$\|u^\lambda_{\bldr\blds}\|=d_\lambda^{-\half}q^{-\psi(\bldr)}$,
% where $\rho$ is the half-sum of positive roots,
% $\lambda(\bldr)$ is the weight such that $e(\lambda,\bldr)$ belongs to the corresponding
% weight space (of $V_\lambda$), and
% $d_\lambda=\sum_{\bldr:\lambda(\bldr)\in P(\lambda)}q^{2(\rho,\lambda(\bldr))}$,
% $P(\lambda)$ being the set of weights corresponding to
% the weight space decomposition of $V_\lambda$.
where
\[
\psi(\bldr)=-\frac{\ell}{2}\sum_{j=1}^{\ell+1}r_{1j}
   + \sum_{i=2}^{\ell+1}\sum_{j=1}^{\ell+2-i}r_{ij},
   \qquad
d_\lambda=\sum_{\bldr:\bldr_1=\lambda} q^{2\psi(\bldr)}
\]
% \psi(\bldr) stands for \pm (\rho,\lambda(\bldr))
%  %    (I'm NOT SURE WHETHER IT IS PLUS OR MINUS!!
% where $\rho$ is the half-sum of positive roots,
% $\lambda(\bldr)$ is the weight such that
% $e(\lambda,\bldr)$ belongs to the corresponding
% weight space (of $V_\lambda$)

Therefore
\be\label{left_mult}
\pi(u_{ij})e^\lambda_{\bldr\blds}
= \sum_{\mu,\bldm,\bldn}
  C_q(\one,\lambda,\mu;i,\bldr,\bldm)C_q(\one,\lambda,\mu;j,\blds,\bldn)
   d_\lambda^\half d_\mu^{-\half}q^{\psi(\bldr)-\psi(\bldm)}
   e^\mu_{\bldm\bldn}.
\ee

Write
\be
\kappa(\bldr,\bldm)=
 d_\lambda^\half d_\mu^{-\half}q^{\psi(\bldr)-\psi(\bldm)}.
\ee
%%%%%%%%%%%%%%%%%%%%%%%%%%%%%%%%%%%%%%%%%%%%%%%%%%%%%%%%%%%%%%%%%%%
\blmma\label{krmbound}
There exist constants $K_2>K_1>0$ such that
$K_1< \kappa(\bldr, M(\bldr))<K_2$ for all $\bldr$.
\elmma
%%%%%%%%%%%%%%%%%%%%%%%%%%%%%%%%%%%%%%%%%%%%%%%%%%%%%%%%%%%%%%%%%%%
\prf
Observe that (\cite{ch-pr}, pp-365)
\[
d_\lambda=\prod_{1\leq i\leq j\leq\ell+1}
\frac{[\lambda_i-\lambda_j+j-i]_q}{[j-i]_q}.
\]
Therefore one gets
\[
\frac{d_\lambda}{d_{\lambda+e_k}}=
  \prod_{j:k<j}\frac{[\lambda_k-\lambda_j+j-k]_q}{[\lambda_k-\lambda_j+j-k+1]_q}
\times
    \prod_{i:i<k}\frac{[\lambda_i-\lambda_k+k-i]_q}{[\lambda_i-\lambda_k+k-i-1]_q}.
\]
There are $\ell$ terms in the above product, and each term
lies between two positive quantities that depend just on $q$.
Next, we have
\[
\psi(\bldr)=-\frac{\ell}{2}\sum_{j=1}^{\ell+1}r_{1j}
   + \sum_{i=2}^{\ell+1}\sum_{j=1}^{\ell+2-i}r_{ij}.
\]
It follows from this that $\psi(\bldr)-\psi(\bldm)$
is bounded.
Therefore the result follows.
\qed

%%%%%%%%%%%%%%%%%%%%%%%%%%%%%%%%%%%%%%%%%%%%%%%%%%%%%%%%%%%%%%
\newsection{Boundedness of commutators}% with equivariant $D$}
%%%%%%%%%%%%%%%%%%%%%%%%%%%%%%%%%%%%%%%%%%%%%%%%%%%%%%%%%%%%%%
Let $D$ be an equivariant Dirac operator acting on $L_2(G)$.
It follows from the discussion in section~3
that $D$ must be of the form
\be
e^\lambda_{\bldr\blds}
\mapsto
d(\bldr)e^\lambda_{\bldr\blds},
\ee
(Here, for a Young tableaux $\lambda$, $N_\lambda$
is the set of all GT tableaux, modulo the appropriate equivalence
relation, with top row $\lambda$).
Then we have
\be\label{bdd_comm}
[D,\pi(u_{ij})]e^\lambda_{\bldr\blds}=
\sum (d(\bldm)-d(\bldr))C_q(\one,\lambda,\mu;i,\bldr,\bldm)
    C_q(\one,\lambda,\mu;j,\blds,\bldn)
\kappa(\bldr,\bldm)e^\mu_{\bldm\bldn}.
\ee
Therefore the condition for boundedness of commutators reads
as follows:
\be \label{eqbdd1}
|(d(\bldm)-d(\bldr))C_q(\one,\lambda,\mu;i,\bldr,\bldm)
   C_q(\one,\lambda,\mu;j,\blds,\bldn)
\kappa(\bldr,\bldm)|<c,
\ee
where $c$ is independent of $i$, $j$, $\lambda$, $\mu$, $\bldr$, $\blds$, $\bldm$ and $\bldn$.

Using lemma~\ref{krmbound}, we get
\be\label{eqbdd2}
|(d(\bldm)-d(\bldr))C_q(\one,\lambda,\mu;i,\bldr,\bldm)
  C_q(\one,\lambda,\mu;j,\blds,\bldn)|<c.
\ee
Choosing $j$, $\blds$ and $\bldn$ suitably, one can ensure that
(\ref{eqbdd2}) implies the following:
\be\label{eqbdd3}
|(d(\bldm)-d(\bldr))C_q(\one,\lambda,\mu;i,\bldr,\bldm)|<c.
\ee
It follows from~(\ref{bdd_comm}) that this condition is also sufficient for
the boundedness of the commutators $[D, u_{ij}]$.

From (\ref{cgc4}), one gets
%%%%%%%%%%%%%%%%%%%%%%%%%%%%%%%%%%%%%%%%%%%%%%%%%%%%%
\be \label{eqbdd4}
|d(\bldr)-d(M(\bldr))|
\leq c q^{-C(i,\bldr,M)}.
\ee
%%%%%%%%%%%%%%%%%%%%%%%%%%%%%%%%%%%%%%%%%%%%%%%%%%%%%

Let us next form a graph $\clg_c$ as described
in section~1 by connecting
two elements $\bldr$ and $\bldr'$ if
$|d(\bldr)-d(\bldr')|<c$.
We will assume
the existence of a partition
$(\Gamma^+,\Gamma^-)$ that does not admit any infinite ladder.
For any subset $F$ of $\Gamma$, we will denote by $F^\pm$ the sets
$F\cap \Gamma^\pm$.
Our next job is to study this graph in more detail
using the boundedness conditions above.
Let us start with  a few definitions and notations.
By an  \textbf{elementary move}, we will mean a map $M$
from some subset of
$\Gamma$ to $\Gamma$ such that $\gamma$ and $M(\gamma)$
are connected by an edge.
A \textbf{move} will mean a composition of a finite number of
elementary moves.
If $M_1$ and $M_2$ are two moves, $M_1M_2$ and $M_2M_1$ will
in general be different.
For a family of moves $M_1, M_2,\ldots, M_r$,
we will denote by
$\sum_{{j=1}}^{r}M_j$
the move $M_1M_2\ldots M_r$,
and by
$\sum_{j=1}^{r}M_{r+1-j}$
the move $M_r\ldots M_2M_1$.
For a nonnegative integer $n$ and a move $M$, we will denote
by $nM$ the move obtained by applying $M$ successively $n$ times.
Of special interest to us will be moves of the
form $M:\bldr\mapsto\blds$, where
$\blds$ is given by (\ref{movenotation}).
We will  use the vector $(m_1,\ldots, m_{k})$
to denote $M$. The following families of moves will
be particularly useful to us:
\[
M_{ik}=(i,i-1,\ldots,i-k+1)\in\bbn^k,\quad
N_{ik}=(\underbrace{i+1,\ldots,i+1}_{\mbox{$k$}},
       i,i,\ldots,i)\in\bbn^{\ell+2-i}.
\]
For describing a path in our graph, we will
often use phrases like `apply the move
$\sum_{{j=1}}^{k}M_j$
to go from $\bldr$ to $\blds$'. This will refer
to the path given by
\[
\Bigl(\bldr,\, M_k(\bldr),
   M_{k-1}M_k(\bldr),\,\ldots,\,M_1M_2\ldots
    M_k(\bldr)=\blds\Bigr).
\]
The following lemma will be very useful in the next
two sections.
%%%%%%%%%%%%%%%%%%%%%%%%%%%%%%%%%%%%%%%%%%%%%%%%%%%%%
\blmma\label{freemove}
Let $N_{jk}$ and $M_{ik}$ be the moves defined above. Then
\begin{enumerate}
\item  $|d(\bldr)-d(N_{j0}(\bldr))|\leq c$,
\item    $|d(\bldr)-d(M_{ik}(\bldr))|\leq
  cq^{-\sum_{a=1}^{k-1}H_{a,i+1-a}-\sum_{b=i}^{\ell}H_{k,b+k-1}}$.  In
  particular, if $H_{a,i+1-a}(\bldr)=0$ for $1\leq a\leq k-1$ and
  $H_{k,b+k-1}(\bldr)=0$ for $i\leq b\leq \ell$, then
  $|d(\bldr)-d(M_{ik}(\bldr))|\leq c$.
\end{enumerate}
\elmma
\prf
Direct consequence of~(\ref{eqbdd4}).
\qed
%%%%%%%%%%%%%%%%%%%%%%%%%%%%%%%%%%%%%%%%%%%%%%%%%%%%%

%%%%%%%%%%%%%%%%%%%%%%%%%%%%%%%%%%%%%%%%%%%%%%%%%%%%%
%%%%%%%%%%%%%%%%%%%%%%%%%%%%%%%%%%%%%%%%%%%%%%%%%%%%%
%%%%%%%%%%%%%%%%%%%%%%%%%%%%%%%%%%%%%%%%%%%%%%%%%%%%%
%%%%%%%%%%%%%%%%%%%%%%%%%%%%%%%%%%%%%%%%%%%%%%%%%%%%%
%%%%%%%%%%%%%%%%%%%%%%%%%%%%%%%%%%%%%%%%%%%%%%%%%%%%%
%%%%%%%%%%%%%%%%%%%%%%%%%%%%%%%%%%%%%%%%%%%%%%%%%%%%%%%%%%%%%%%%%
\newsection{Characterization of $|D|$}
%%%%%%%%%%%%%%%%%%%%%%%%%%%%%%%%%%%%%%%%%%%%%%%%%%%%%%%%%%%%%%%%%
In this section and the next, we will use lemma~\ref{freemove}
to prove a characterization theorem for the sign of the
operator $D$. Along the way,
we will also give a very precise description
of the singular values of $D$.
The main ingredients in the proof are
the finiteness of exactly one of the sets $F^+$ and $F^-$
for appropriately chosen subsets $F$ of $\Gamma$.
General form of the argument for proving this will be
as follows:
for a carefully chosen coordinate $C$ (in the present case, $C$
would be one of the $V_{a1}$'s or $H_{ab}$'s), a sweepout argument
will show that any $\gamma$ can be connected by a path,
throughout which $C(\cdot)$ remains constant, to another point $\gamma'$
for which $C(\gamma')=C(\gamma)$ and all other coordinates of $\gamma'$
are zero.
This would help connect any two points $\gamma$ and $\delta$ by a path
such that $C(\cdot)$ would lie between $C(\gamma)$ and $C(\delta)$
on the path. This would finally result in the finiteness
of at least one (and hence exactly one) of $C(F^+)$ and $C(F^-)$.
Next, assuming one of these, say $C(F^-)$ is finite,
one shows that for any other coordinate $C'$,
$C'(F^-)$ is also finite.
This is done as follows. If $C'(F^-)$ is infinite, one chooses
elements $y_n\in F^-$ with
$C'(y_n)<C'(y_{n+1})$ for all $n$.
Now starting at each $y_n$, produce paths
keeping the $C'$-coordinate constant and taking the
$C$-coordinate above the plane $C(\cdot)=K$, where $C(F^-)\seq [-K,K]$.
This will produce an infinite ladder.
The argument is explained in the following diagram.\\[3ex]
%%%%%%%%%%%%%%%%%%%%%%%%%%%%%%%%%%%%%%%%%%%%%%%%%%%%
\hspace*{60pt}
%%%%%%%%%%%%%%%%%%%%%%%%%%%%%%%%%%%%%%%%%%%%%%%%%%%%
\setlength{\unitlength}{0.00041667in}
\begingroup\makeatletter\ifx\SetFigFont\undefined%
\gdef\SetFigFont#1#2#3#4#5{%
  \reset@font\fontsize{#1}{#2pt}%
  \fontfamily{#3}\fontseries{#4}\fontshape{#5}%
  \selectfont}%
\fi\endgroup%
{\renewcommand{\dashlinestretch}{30}
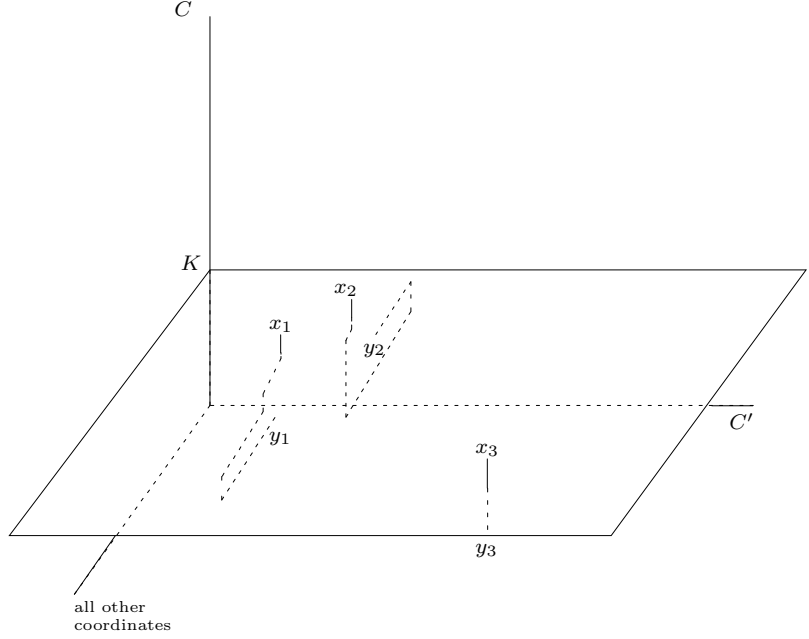
\begin{picture}(10149,7971)(0,-10)
\path(2562,2850)(2562,7800)
\dashline{60.000}(2562,2850)(837,450)
\dashline{60.000}(9462,2850)(2562,2850)(2562,4575)
\dashline{60.000}(4287,3675)(4287,2700)(5112,4050)
        (5112,4425)(4662,3750)
\path(1362,1200)(837,450)
\path(8937,2850)(9462,2850)
\path(2562,4575)(12,1200)(7662,1200)
        (10137,4575)(2562,4575)
\path(3462,3750)(3462,3525)
\dashline{60.000}(3462,3525)(3462,3450)(3237,3000)
        (3237,2775)(2712,1950)(2712,1650)(3387,2700)
\dashline{60.000}(4287,3675)(4362,3825)(4362,3975)
\path(4362,3975)(4362,4200)
\dashline{60.000}(6087,1275)(6087,1800)
\path(6087,1800)(6087,2175)
\put(2112,7800){\makebox(0,0)[lb]{\smash{{{\SetFigFont{8}{9.6}{\rmdefault}{\mddefault}{\updefault}$C$}}}}}
\put(2187,4575){\makebox(0,0)[lb]{\smash{{{\SetFigFont{8}{9.6}{\rmdefault}{\mddefault}{\updefault}$K$}}}}}
\put(837,225){\makebox(0,0)[lb]{\smash{{{\SetFigFont{6}{7.2}{\rmdefault}{\mddefault}{\updefault}all other}}}}}
\put(837,0){\makebox(0,0)[lb]{\smash{{{\SetFigFont{6}{7.2}{\rmdefault}{\mddefault}{\updefault}coordinates}}}}}
\put(9162,2550){\makebox(0,0)[lb]{\smash{{{\SetFigFont{8}{9.6}{\rmdefault}{\mddefault}{\updefault}$C'$}}}}}
\put(4512,3525){\makebox(0,0)[lb]{\smash{{{\SetFigFont{8}{9.6}{\rmdefault}{\mddefault}{\updefault}$y_2$}}}}}
\put(3312,2400){\makebox(0,0)[lb]{\smash{{{\SetFigFont{8}{9.6}{\rmdefault}{\mddefault}{\updefault}$y_1$}}}}}
\put(3312,3825){\makebox(0,0)[lb]{\smash{{{\SetFigFont{8}{9.6}{\rmdefault}{\mddefault}{\updefault}$x_1$}}}}}
\put(5937,975){\makebox(0,0)[lb]{\smash{{{\SetFigFont{8}{9.6}{\rmdefault}{\mddefault}{\updefault}$y_3$}}}}}
\put(5937,2250){\makebox(0,0)[lb]{\smash{{{\SetFigFont{8}{9.6}{\rmdefault}{\mddefault}{\updefault}$x_3$}}}}}
\put(4137,4275){\makebox(0,0)[lb]{\smash{{{\SetFigFont{8}{9.6}{\rmdefault}{\mddefault}{\updefault}$x_2$}}}}}
\end{picture}
}\\[2ex]
%%%%%%%%%%%%%%%%%%%%%%%%%%%%%%%%%%%%%%%%%%%%%%%%%%%%%%%%%%%%%%%%

Our next job is to define an important class of subsets of $\Gamma$.
Observe that lemma~\ref{freemove}
tells us that for any $\bldr$ and any $j$, the points
$\bldr$ and $N_{j0}(\bldr)$ are connected by an edge,
whenever $N_{j0}(\bldr)$ is a GT tableaux.
Let $\bldr$ be an element of $\Gamma$. Define the
\textbf{free plane passing through $\bldr$} to be  the minimal
subset of $\Gamma$ that contains $\bldr$
and is closed under application of the moves $N_{j0}$.
We will denote this set by $\scrf_\bldr$.
The following is an easy consequence of this definition.
%%%%%%%%%%%%%%%%%%%%%%%%%%%%%%%%%%%%%%%%%%%%%%%%%
\blmma \label{freecriterion}
Let $\bldr$ and $\blds$ be two GT tableaux. Then
$\blds\in \scrf_\bldr$ if and only if
$V_{a,1}(\bldr)=V_{a,1}(\blds)$ for all $a$ and for each $b$, the difference
$H_{a,b}(\bldr)-H_{a,b}(\blds)$ is independent of $a$.
\elmma
%%%%%%%%%%%%%%%%%%%%%%%%%%%%%%%%%%%%%%%%%%%%%%%%%
%%%%%%%%%%%%%%%%%%%%%%%%%%%%%%%%%%%%%%%%%%%%%%%%%
\bcrlre \label{freedisjt}
Let $\bldr,\blds\in\Gamma$. Then either $\scrf_\bldr=\scrf_\blds$ or
$\scrf_\bldr\cap \scrf_\blds=\phi$.
\ecrlre
%%%%%%%%%%%%%%%%%%%%%%%%%%%%%%%%%%%%%%%%%%%%%%%%%

Let $\bldr\in\Gamma$.
For $1\leq j\leq \ell+1$, define
$a_j$ to be an integer such that $H_{a_j,j}(\bldr)=\min_i H_{ij}(\bldr)$.
Note three things here:\\
1. definition of $a_j$ depends on $\bldr$,\\
2. for a given $j$ and given $\bldr$, $a_j$ need not be unique, and\\
3. if $\blds\in\scrf_\bldr$, then for each $j$, the set of
$k$'s for which $H_{kj}(\blds)=\min_i H_{ij}(\blds)$ is same
as the set of all $k$'s for which $H_{kj}(\bldr)=\min_i H_{ij}(\bldr)$.
Therefore, the $a_j$'s can be chosen in a manner such that
they remain the same for all elements lying on a given free plane.

%%%%%%%%%%%%%%%%%%%%%%%%%%%%%%%%%%%%%%%%%%%%%%%%%%%%
\blmma\label{sweep1}
Let $\blds\in \scrf_\bldr$. Let $\blds'$ be another GT tableaux
given by
\[
V_{a1}(\blds')=V_{a1}(\blds) \mbox{ and }
H_{a1}(\blds')=H_{a1}(\blds) \mbox{ for all }a,\quad
 H_{a_b,b}(\blds')=0 \mbox{ for all }b>1,
\]
where the $a_j$'s are as defined above.
Then there is a path in $\scrf_\bldr$ from $\blds$ to $\blds'$
such that $H_{11}(\cdot)$ remains constant throughout this path.
\elmma
%%%%%%%%%%%%%%%%%%%%%%%%%%%%%%%%%%%%%%%%%%%%%%%%%%%%
\prf
Apply the move
$\sum_{{b=2}}^{\ell}
   \left(\sum_{j=2}^{\ell+2-b} H_{a_j,j}(\blds)\right)N_{\ell+3-b,0}$.\qed

The following diagram will help explain the steps involved
in the above proof in the case where $\bldr$ is the constant
tableaux.\\[2ex]
\def\labelstyle{\scriptstyle}
\xymatrix@C=.6pt@R=.6pt{
    \cdot\ar@{}[r] & & \cdot\ar@{}[r] && \odot\ar@{.}[d] && \cdot&& \cdot&\\
    0 & a &  & b & & c && d &\\
    \cdot\ar@{}[r] & & \cdot\ar@{}[r] && \cdot\ar@{.}[u]\ar@{.}[d] && \cdot&\\
    0 & a &  & b & & c &\\
    \cdot\ar@{}[r] & & \cdot\ar@{}[r] && \odot\ar@{.}[u]&\\
    0 & a && b &  \\
    \cdot\ar@{}[r] && \cdot&\\
    0 &a & \\
    \cdot&}
\hspace{-2em} \xymatrix@C=20pt@R=12pt{&\\&\\ \ar@{->}[r]^{bN_{30}}&\\}\hspace{.3em}
\xymatrix@C=.6pt@R=.6pt{
    \cdot\ar@{}[r] & & \cdot\ar@{}[r] && \cdot && \odot\ar@{.}[d]&& \cdot&\\
    0 & a &  & 0 & & b+c && d &\\
    \cdot\ar@{}[r] & & \cdot\ar@{}[r] && \cdot && \odot\ar@{.}[u]&\\
    0 & a &  & 0 & & b+c &\\
    \cdot\ar@{}[r] & & \cdot\ar@{}[r] && \cdot&\\
    0 & a && 0 &  \\
    \cdot\ar@{}[r] && \cdot&\\
    0 &a & \\
    \cdot&}
\hspace{-2em} \xymatrix@C=20pt@R=12pt{&\\&\\ \ar@{->}[r]^{(b+c)N_{40}}&\\}\hspace{.3em}
\xymatrix@C=.6pt@R=.6pt{
    \cdot\ar@{}[r] & & \cdot\ar@{}[r] && \cdot && \cdot& \odot&\\
    0 & a &  & 0 & & 0 && b+c+d &\\
    \cdot\ar@{}[r] & & \cdot\ar@{}[r] && \cdot && \cdot&\\
    0 & a &  & 0 & & 0 &\\
    \cdot\ar@{}[r] & & \cdot\ar@{}[r] && \cdot&\\
    0 & a && 0 &  \\
    \cdot\ar@{}[r] && \cdot&\\
    0 &a & \\
    \cdot&}\\[2ex]
  \hspace*{12em}  \xymatrix@C=20pt@R=12pt{&\\&\\ \ar@{->}[r]^{(b+c+d)N_{50}}&\\}\hspace{.3em}
\xymatrix@C=.6pt@R=.6pt{
    \cdot\ar@{}[r] & & \cdot\ar@{}[r] && \cdot && \cdot&& \cdot&\\
    0 & a &  & 0 & & 0 && 0 &\\
    \cdot\ar@{}[r] & & \cdot\ar@{}[r] && \cdot && \cdot&\\
    0 & a &  & 0 & & 0 &\\
    \cdot\ar@{}[r] & & \cdot\ar@{}[r] && \cdot&\\
    0 & a && 0 &  \\
    \cdot\ar@{}[r] && \cdot&\\
    0 &a & \\
    \cdot&}\\
A dotted line joining two circled dots signifies a move that
increases the $r_{ij}$'s lying on the dotted line by one.
Where there is one circled dot and no dotted line, it means
one applies the move that raises the $r_{ij}$ corresponding to
the circled dot by one.

%%%%%%%%%%%%%%%%%%%%%%%%%%%%%%%%%%%%%%%%%%%%%%%%%%%%
\bppsn\label{signfree1}
Let $\bldr$ be a GT tableaux.
Then either $\scrf_{\bldr}^+$ is finite or $\scrf_{\bldr}^-$ is finite.
\eppsn
%%%%%%%%%%%%%%%%%%%%%%%%%%%%%%%%%%%%%%%%%%%%%%%%%%%%
\prf
Suppose, if possible, both $H_{11}(\scrf_{\bldr}^+)$ and $H_{11}(\scrf_{\bldr}^-)$
are infinite. Then there exist two sequences of elements
$\bldr_n$ and $\blds_n$ with $\bldr_n\in \scrf_\bldr^+$
and $\blds_n\in \scrf_\bldr^-$,
such that
\[
H_{11}(\bldr_1)<H_{11}(\blds_1)<H_{11}(\bldr_2)<H_{11}(\blds_2)<\cdots.
\]
Now starting from $\bldr_n$,  employ the forgoing lemma
to reach a point $\bldr'_n\in\scrf_{\bldr}$ for which
\[
V_{a1}(\bldr'_n)=V_{a1}(\bldr_n) \mbox{ and }
H_{a1}(\bldr'_n)=H_{a1}(\bldr_n) \mbox{ for all }a,\quad
 H_{a_b,b}(\bldr'_n)=0 \mbox{ for all }b>1.
\]
Similarly, start at $\blds_n$ and go to
a point $\blds'_n\in\scrf_{\bldr}$ for which
\[
V_{a1}(\blds'_n)=V_{a1}(\blds_n) \mbox{ and }
H_{a1}(\blds'_n)=H_{a1}(\blds_n) \mbox{ for all }a,\quad
 H_{a_b,b}(\blds'_n)=0 \mbox{ for all }b>1.
\]
Now use the move $N_{10}$ to get to $\blds'_n$ from $\bldr'_n$.
The paths thus constructed are all disjoint, because
for the path from $\bldr_n$ to $\blds_n$, the
$H_{11}$ coordinate lies between
$H_{11}(\bldr_n)$ and $H_{11}(\blds_n)$.
This means $(\scrf_\bldr^+, \scrf_\bldr^-)$ admits an infinite ladder.
So one of the sets  $H_{11}(\scrf_\bldr^+)$ and $H_{11}(\scrf_\bldr^-)$
must be finite. Let us assume that $H_{11}(\scrf_\bldr^-)$ is finite.

Let us next show that for any $b>1$, $H_{ab}(\scrf_\bldr^-)$ is finite.
Let $K$ be an integer such that $H_{11}(\blds)<K$ for all $\blds\in \scrf_\bldr^-$.
If $H_{ab}(\scrf_\bldr^-)$ was infinite, there would exist elements
$\bldr_n\in \scrf_\bldr^-$ such that
\[
H_{ab}(\bldr_1)<H_{ab}(\bldr_2)<\cdots.
\]
Now start at $\bldr_n$ and employ the move $N_{10}$ successively
$K$ times to reach a point in
$\scrf_\bldr^+=\scrf_\bldr\backslash\scrf_\bldr^-$.
These paths will all be disjoint, as throughout the path,
$H_{ab}$ remains fixed.

Since the coordinates $(H_{11},H_{12},\ldots,H_{1,\ell})$
completely specify a point in $\scrf_\bldr$, it follows that
$\scrf_\bldr^-$ is finite.
\qed

Next we need a set that can be used for a proper
indexing of the free planes.
Such a set will be called a complementary axis.

%%%%%%%%%%%%%%%%%%%%%%%%%%%%%%%%%%%%%%%%%%%%%%%%%%%%
\bdfn\rm
A subset $\scrc $ of $\Gamma$ is called a \textbf{complementary axis} if
\begin{enumerate}
\item $\cup_{\bldr\in \scrc }\scrf_\bldr =\Gamma$,
\item if $\bldr,\blds\in \scrc $, and $\bldr\neq \blds$, then
   $\scrf_\bldr$ and $\scrf_\blds$ are disjoint.
\end{enumerate}
\edfn
%%%%%%%%%%%%%%%%%%%%%%%%%%%%%%%%%%%%%%%%%%%%%%%%%%%%

Let us next give a choice of a complementary axis.
%%%%%%%%%%%%%%%%%%%%%%%%%%%%%%%%%%%%%%%%%%%%%%%%%%%%
\bthm \label{compl}
Define
\[
\scrc =\{\bldr\in \Gamma: \Pi_{a=1}^{\ell+1-b} H_{ab}(\bldr)=0
           \mbox{ for } 1\leq b\leq \ell\}.
\]
The set $\scrc $ defined above is a complementary axis.
\ethm
%%%%%%%%%%%%%%%%%%%%%%%%%%%%%%%%%%%%%%%%%%%%%%%%%%%%
\prf
Let $\blds\in\Gamma$.
A sweepout argument almost identical to that used in
lemma~\ref{sweep1} (application of the move
$\sum_{{b=1}}^\ell
   \left(\sum_{j=1}^{\ell+1-b} H_{a_j,j}(\blds)\right)N_{\ell+2-b,0}$ )
will connect $\blds$ to another element $\blds'$
for which $H_{a_b,b}(\blds')=0$ for $1\leq b\leq\ell$
by a path that lies entirely on $\scrf_\blds$.
Clearly, $\blds'\in\scrc$. Since $\blds'\in\scrf_\blds$,
by corollary~\ref{freedisjt}, $\blds\in\scrf_{\blds'}$.

It remains to show that if $\bldr$ and $\blds$ are two distinct elements of
$\scrc$, then $\blds\not\in\scrf_\bldr$.
Since $\bldr\neq\blds$, there exist two integers  $a$ and $b$,
$1\leq b\leq \ell$ and $1\leq a\leq \ell+2-b$, such that
$H_{ab}(\bldr)\neq H_{ab}(\blds)$.
Observe that $H_{1\ell}(\cdot)$ must be zero for both,
as they are members of $\scrc$. So $b$ can not be $\ell$ here.
Next we will produce two integers $i$ and $j$ such that
the differences $H_{ib}(\bldr)-H_{ib}(\blds)$
and $H_{jb}(\bldr)-H_{jb}(\blds)$ are distinct.
If there is an integer $k$ for which
$H_{kb}(\bldr)=H_{kb}(\blds)=0$, then take $i=a$, $j=k$.
If not, there would exist two integers $i$ and $j$ such that
$H_{ib}(\bldr)=0$, $H_{ib}(\blds)>0$
and
$H_{jb}(\bldr)>0$, $H_{jb}(\blds)=0$.
Take these $i$ and $j$.
Since $H_{ib}(\bldr)-H_{ib}(\blds)$
and $H_{jb}(\bldr)-H_{jb}(\blds)$ are distinct,
by lemma~\ref{freecriterion}, $\bldr$ and $\blds$ can not lie
on the same free plane.
\qed

%%%%%%%%%%%%%%%%%%%%%%%%%%%%%%%%%%%%%%%
\blmma \label{sweep2}
Let $\bldr$ be a GT tableaux. Let $\blds$ be the GT tableaux
defined by the prescription
\[
V_{a1}(\blds)=V_{a1}(\bldr)\mbox{ for all }a,\quad
 H_{ab}(\blds)=H_{ab}(\bldr) \mbox{ for all }a\geq 2,\mbox{ for all }b,\quad
 H_{1,b}(\blds)=0 \mbox{ for all }b.
\]
Then there is a path from $\bldr$ to $\blds$ such that
$V_{a1}(\cdot)$ remains constant throughout the path.
\elmma
%%%%%%%%%%%%%%%%%%%%%%%%%%%%%%%%%%%%%%%
\prf
Apply  the move
$\displaystyle{\sum_{{b=1}}^\ell} H_{1,b}(\bldr)M_{b+1,1}$.\qed

The above lemma is actually the first step in the following
slightly more general sweepout algorithm.

%%%%%%%%%%%%%%%%%%%%%%%%%%%%%%%%%%%%%%%%%%%%%%%%%%%%
\blmma \label{sweep3}
Let $\bldr$ be a GT tableaux. Let $\blds$ be the GT tableaux
defined by the prescription
\[
V_{11}(\blds)=V_{11}(\bldr),\quad
V_{a1}(\blds)=0\mbox{ for all }a>1,\quad
 H_{ab}(\blds)=0 \mbox{ for all }a,b.
\]
Then there is a path from $\bldr$ to $\blds$ such that
$V_{11}(\cdot)$ remains constant throughout the path.
\elmma
%%%%%%%%%%%%%%%%%%%%%%%%%%%%%%%%%%%%%%%%%%%%%%%%%%%%
\prf
Apply successively the moves
\[
\sum_{{b=1}}^\ell H_{1,b}(\bldr)M_{b+1,1},\quad
\sum_{{b=1}}^{\ell-1} H_{2,b}(\bldr)M_{b+2,2},\quad
\ldots,\quad
 H_{\ell,1}(\bldr)M_{\ell+1,\ell},
\]
followed by
%%%%%%%%%%%%%%%%%%%%%%%%%%%%%%%%%%%%%%%%%%%%%%%%%%%%
\be\label{movseq}
V_{21}(\bldr)M_{33},\quad  (V_{21}(\bldr)+V_{31}(\bldr))M_{44},\quad
\ldots,\quad  \left(\sum_{a=2}^\ell V_{a1}(\bldr)\right)M_{\ell+1,\ell+1}.
\ee
%%%%%%%%%%%%%%%%%%%%%%%%%%%%%%%%%%%%%%%%%%%%%%%%%%%%
\qed

The following diagram
will help explain the procedure described above in a simple case.\\
%%%%%%%%%%%%%%%%%%%%%%%%%%%%%%%%%%%%%
\def\labelstyle{\scriptstyle}
\xymatrix@C=10pt@R=8pt{
    \cdot\ar@{}[r]\ar@{}[d]|\star      & \cdot\ar@{}[r] & \cdot\ar@{}[r] & \odot&\\
    \cdot\ar@{}[d]|\star\ar@{}[ur]|\star& \cdot\ar@{}[ur]|\star & \cdot\ar@{}[ur]|\star &  \\
    \cdot\ar@{}[d]|\star\ar@{}[ur]|\star& \cdot\ar@{}[ur]|\star & \\
     \cdot\ar@{}[ur]|\star& }
 \hspace{-1em} \xymatrix@C=20pt@R=12pt{&\\ \ar@{->}[r]^{M_{41}}&\\}\hspace{.5em}
\xymatrix@C=10pt@R=8pt{
    \cdot\ar@{}[r]\ar@{}[d]|\star      & \cdot\ar@{}[r] & \odot\ar@{}[r] & \cdot&\\
    \cdot\ar@{}[d]|\star\ar@{}[ur]|\star& \cdot\ar@{}[ur]|\star & \cdot\ar@{}[ur]|0 &  \\
    \cdot\ar@{}[d]|\star\ar@{}[ur]|\star& \cdot\ar@{}[ur]|\star & \\
     \cdot\ar@{}[ur]|\star&  }
  \hspace{-1em}  \xymatrix@C=20pt@R=12pt{&\\ \ar@{->}[r]^{M_{31}}&\\}\hspace{.5em}
\xymatrix@C=10pt@R=8pt{
    \cdot\ar@{}[r]\ar@{}[d]|\star      & \odot\ar@{}[r] & \cdot\ar@{}[r] & \cdot&\\
    \cdot\ar@{}[d]|\star\ar@{}[ur]|\star& \cdot\ar@{}[ur]|0 & \cdot\ar@{}[ur]|0 &  \\
    \cdot\ar@{}[d]|\star\ar@{}[ur]|\star& \cdot\ar@{}[ur]|\star & \\
     \cdot\ar@{}[ur]|\star& \\
     }
  \hspace{-1em}  \xymatrix@C=20pt@R=12pt{&\\ \ar@{->}[r]^{M_{21}}&\\}\hspace{.5em}
\xymatrix@C=10pt@R=8pt{
    \cdot\ar@{}[r]\ar@{}[d]|\star      & \cdot\ar@{}[r] & \cdot\ar@{}[r] & \odot&\\
    \cdot\ar@{}[d]|\star\ar@{}[ur]|0& \cdot\ar@{}[ur]|0 & \odot\ar@{.}[ur]|0 &  \\
    \cdot\ar@{}[d]|\star\ar@{}[ur]|\star& \cdot\ar@{}[ur]|\star & \\
     \cdot\ar@{}[ur]|\star& \\
     }\\[2ex]
\hspace*{60pt}
  \xymatrix@C=20pt@R=12pt{&\\ \ar@{->}[r]^{M_{42}}&\\}\hspace{.5em}
\xymatrix@C=10pt@R=8pt{
    \cdot\ar@{}[r]\ar@{}[d]|\star      & \cdot\ar@{}[r] & \odot\ar@{}[r] & \cdot&\\
    \cdot\ar@{}[d]|\star\ar@{}[ur]|0& \odot\ar@{.}[ur]|0 & \cdot\ar@{}[ur]|0 &  \\
    \cdot\ar@{}[d]|\star\ar@{}[ur]|\star& \cdot\ar@{}[ur]|0 & \\
     \cdot\ar@{}[ur]|\star& \\
     }
  \hspace{-1em}  \xymatrix@C=20pt@R=12pt{&\\ \ar@{->}[r]^{M_{32}}&\\}\hspace{.5em}
\xymatrix@C=10pt@R=8pt{
    \cdot\ar@{}[r]\ar@{}[d]|\star      & \cdot\ar@{}[r] & \cdot\ar@{}[r] & \odot&\\
    \cdot\ar@{}[d]|\star\ar@{}[ur]|0& \cdot\ar@{}[ur]|0 & \cdot\ar@{.}[ur]|0 &&  \\
    \cdot\ar@{}[d]|\star\ar@{}[ur]|0& \odot\ar@{.}[ur]|0 &&& \\
     \cdot\ar@{}[ur]|\star& &&&\\
     }
  \hspace{-1em}  \xymatrix@C=20pt@R=12pt{&\\ \ar@{->}[r]^{M_{43}}&\\}\hspace{.5em}
\xymatrix@C=10pt@R=8pt{
    \cdot\ar@{}[r]\ar@{}[d]|\star      & \cdot\ar@{}[r] & \odot\ar@{}[r] & \cdot&\\
    \cdot\ar@{}[d]|\star\ar@{}[ur]|0& \cdot\ar@{.}[ur]|0 & \cdot\ar@{}[ur]|0 &  \\
    \odot\ar@{}[d]|\star\ar@{.}[ur]|0& \cdot\ar@{}[ur]|0 & \\
     \cdot\ar@{}[ur]|0& \\
     }\\[2ex]
\hspace*{60pt}
  \xymatrix@C=20pt@R=12pt{&\\ \ar@{->}[r]^{M_{33}}&\\}\hspace{.5em}
\xymatrix@C=10pt@R=8pt{
    \cdot\ar@{}[r]\ar@{}[d]|\star      & \cdot\ar@{}[r] & \cdot\ar@{}[r] & \odot&\\
    \cdot\ar@{}[d]|0\ar@{}[ur]|0& \cdot\ar@{}[ur]|0 & \cdot\ar@{.}[ur]|0 &  \\
    \cdot\ar@{}[d]|\star\ar@{}[ur]|0& \cdot\ar@{.}[ur]|0 & \\
     \odot\ar@{.}[ur]|0& \\
     }
  \xymatrix@C=20pt@R=12pt{&\\ \ar@{->}[r]^{M_{44}}&\\}\hspace{.5em}
\xymatrix@C=10pt@R=8pt{
    \cdot\ar@{}[r]\ar@{}[d]|\star      & \cdot\ar@{}[r] & \cdot\ar@{}[r] & \cdot&\\
    \cdot\ar@{}[d]|0\ar@{}[ur]|0& \cdot\ar@{}[ur]|0 & \cdot\ar@{}[ur]|0 &  \\
    \cdot\ar@{}[d]|0\ar@{}[ur]|0& \cdot\ar@{}[ur]|0 & \\
     \cdot\ar@{}[ur]|0& \\
     }

%%%%%%%%%%%%%%%%%%%%%%%%%%%%%%%%%%%%%%%%%
\bcrlre\label{growth5}
$|d(\bldr)|=O(r_{11})$.
\ecrlre
%%%%%%%%%%%%%%%%%%%%%%%%%%%%%%%%%%%%%%%%%
\prf
If one employs the sequence of moves
\[
V_{11}(\bldr)M_{22},\quad (V_{11}(\bldr)+V_{21}(\bldr))M_{33},\quad
\ldots,\quad  \left(\sum_{a=1}^\ell V_{a1}(\bldr)\right)M_{\ell+1,\ell+1}
\]
instead of the sequence given in (\ref{movseq}),
one would reach the constant (or zero) tableaux.
Total length of this path from $\bldr$ to the zero tableaux
is
\[
\sum_{a=1}^\ell\sum_{b=1}^{\ell+1-a}H_{ab}(\bldr)
   + \sum_{b=1}^\ell\sum_{a=1}^b V_{a1}(\bldr),
\]
which can easily be shown to be bounded by $\ell r_{11}$.
\qed

%%%%%%%%%%%%%%%%%%%%%%%%%%%%%%%%%%%%%%%%%
\bthm\label{singular}
Let $\widetilde{D}$  be the following operator:
\be
\widetilde{D}: e^\lambda_{\bldr,\blds}\mapsto r_{11}e^\lambda_{\bldr,\blds}
\ee
Then $(\cla,\clh,\widetilde{D})$ is an equivariant $\ell(\ell+2)$-summable odd spectral triple.

Moreover, if $D$ is any equivariant Dirac operator
acting on the $L_2$ space of $SU_q(\ell+1)$, then
there exist positive reals $a$ and $b$ such that
$|D| \leq a + b\widetilde{D}$.
In particular,
$D$ cannot be $p$-summable for $p<\ell(\ell+2)$.
\ethm
%%%%%%%%%%%%%%%%%%%%%%%%%%%%%%%%%%%%%%%%%
\prf
%%%%%%%%%%%%%%%%%%%%%
Boundedness of commutators with algebra elements
follow from the observation that
$|d(\bldr)-d(M(\bldr)|\leq 1$ and hence
equation~(\ref{eqbdd3}) is satisfied.

Observe that the number of Young tableux
$\lambda=(\lambda_1,\ldots,\lambda_\ell,\lambda_{\ell+1})$
with
$n=\lambda_1\geq \lambda_2\geq \ldots \lambda_\ell\geq \lambda_{\ell+1}=0$
 is
\[
\sum_{i_1=0}^{n}\sum_{i_2=0}^{i_1}\ldots\sum_{i_{\ell-1}=0}^{i_{\ell-2}}1
=
\mbox{polynomial in $n$ of degree $\ell-1$}.
\]
Thus the number of such Young tableaux is $O(n^{\ell-1})$.

Next, let
$\lambda:n=\lambda_1\geq \lambda_2\geq\ldots\geq \lambda_{\ell}\geq 0$
be an Young tableaux, and let $V_\lambda$ be the space
carrying the irreducible representation parametrized by
$\lambda$. Then
\bean
\mbox{dim}\, V_\lambda &=&
 \prod_{1\leq i<j\leq \ell+1}
  \frac{(\lambda_{i}-\lambda_{i+1})+
      \ldots (\lambda_{j-1}-\lambda_{j})+j-i}{j-i}\\
  &=& \prod_{1\leq i<j\leq \ell+1}
  \frac{\lambda_{i}-\lambda_{j}+j-i}{j-i}\\
  &\leq& (n+1)^{\frac{\ell(\ell+1)}{2}}.
\eean
Thus the dimension of an irreducible representation corresponding to
a Young tableaux
\[
n=\lambda_1\geq \lambda_2\geq \ldots \lambda_\ell\geq \lambda_{\ell+1}=0
\]
is $O(n^{\half\ell(\ell+1)})$.

Using the two observations above, one can now show
that the summability
of $\widetilde{D}$ is $\ell(\ell+2)$.
Optimality of $\widetilde{D}$ follows from
corollary~\ref{growth5}.\qed

One should note, however, that the $\widetilde{D}$ defined
above has trivial sign, and consequently trivial $K$-homology class.

%% optimality of \widetilde{D}\otimes I
\blmma\label{optimality1}
Let $D$ be an equivariant Dirac operator
on $L_2(G)\otimes\bbc^m$. Then
there are positive reals $a,b$  such that
$|D|\leq a+b|\widetilde{D}\otimes I|$.
\elmma
\prf
Let $D$ be an equivariant Dirac operator on $L_2(G)\otimes\bbc^m$.
Then $D$ must be of the form
$e_{\bldr,\blds}\otimes v\mapsto e_{\bldr,\blds}\otimes T(\bldr)v$
where $T(\bldr)$ are self-adjoint operators acting on $\bbc^m$.
The growth conditions coming out of the boundedness of
the commutators will now be exactly as in~(\ref{eqbdd4}), with
the scalars $d(\cdot)$ replaced by operators $T(\cdot)$ and
absolute value replaced by operator norm.
If we now form a  graph by joining two vertices
$\bldr$ and $\blds$ whenever $\|T(\bldr)-T(\blds)\|\leq c$,
then exactly as in the proof of corollary~\ref{growth5},
one can show that any point $\bldr$ can be connected to
the zero tableaux by a path of length $O(r_{11})$.
This implies that there are positive reals $a$ and $b$ such
that $|T(\bldr)|\leq a+br_{11}$.
The assertion in the lemma now follows from this.\qed
%%%%%%%%%%%%%%%%%%%%%%%%%%%%%%%%%%%%%%%%%%%%%%%%%%%%%
%%%%%%%%%%%%%%%%%%%%%%%%%%%%%%%%%%%%%%%%%%%%%%%%%%%%%
%%%%%%%%%%%%%%%%%%%%%%%%%%%%%%%%%%%%%%%%%%%%%%%%%%%%%
%%%%%%%%%%%%%%%%%%%%%%%%%%%%%%%%%%%%%%%%%%%%%%%%%%%%%
%%%%%%%%%%%%%%%%%%%%%%%%%%%%%%%%%%%%%%%%%%%%%%%%%%%%%
%%%%%%%%%%%%%%%%%%%%%%%%%%%%%%%%%%%%%%%%%%%%%%%%%%%%%
%%%%%%%%%%%%%%%%%%%%%%%%%%%%%%%%%%%%%%%%%%%%%%%%%%%%%
%%%%%%%%%%%%%%%%%%%%%%%%%%%%%%%%%%%%%%%%%%%%%%%%%%%%%
%%%%%%%%%%%%%%%%%%%%%%%%%%%%%%%%%%%%%%%%%%%%%%%%%%%%%
%%%%%%%%%%%%%%%%%%%%%%%%%%%%%%%%%%%%%%%%%%%%%%%%%%%%%
%%%%%%%%%%%%%%%%%%%%%%%%%%%%%%%%%%%%%%%%%%%%%%%%%%%%%
%%%%%%%%%%%%%%%%%%%%%%%%%%%%%%%%%%%%%%%%%
\newsection{Characterization of $\mbox{sign}\,D$}
%%%%%%%%%%%%%%%%%%%%%%%%%%%%%%%%%%%%%%%%%%%%%%%%%%%%%
We continue our analysis of the growth
conditions on the $d(\bldr)$'s in this section
in order to come up with a complete characterization
of the sign of $D$.
%%%%%%%%%%%%%%%%%%%%%%%%%%%%%%%%%%%%%%%%%
\blmma\label{about_v11}
The sets $V_{11}(\Gamma^+)$ and $V_{11}(\Gamma^-)$ can not both be infinite.
\elmma
%%%%%%%%%%%%%%%%%%%%%%%%%%%%%%%%%%%%%%%%%
\prf
If both the sets are infinite, then one can
choose two sequences of points $\bldr_n$ and $\blds_n$
such that $\bldr_n\in \Gamma^+$, $\blds_n\in \Gamma^-$ and
\[
V_{11}(\bldr_1)<V_{11}(\blds_1)<  V_{11}(\bldr_2)<  V_{11}(\blds_2)<\ldots.
\]
Start at $\bldr_n$ and use lemma~\ref{sweep3} above to reach a point
$\bldr'_n$ for which $V_{11}(\bldr'_n)=V_{11}(\bldr_n)$ and all other coordinates
are zero through a path where the $V_{11}$ coordinate remains constant.
Similarly, from $\blds_n$, go to a point $\blds'_n$ for which
$V_{11}(\blds'_n)=V_{11}(\blds_n)$ and all other coordinates
are zero. Now apply  the move
$(V_{11}(\blds_n)-V_{11}(\bldr_n))M_{11}$ to go from $\bldr'_n$ to $\blds'_n$.
This will give us a path $p_n$ from $\bldr_n$ to $\blds_n$
on which $V_{11}(\cdot)$ remains between $V_{11}(\bldr_n)$
and $V_{11}(\blds_n)$. Therefore all the paths $p_n$ are disjoint.
Thus $(\Gamma^+,\Gamma^-)$ admits an infinite ladder.
So at least one of $V_{11}(\Gamma^+)$ and $V_{11}(\Gamma^-)$
must be finite.\qed

%%%%%%%%%%%%%%%%%%%%%%%%%%%%%%%%%%%%%%%%%
\blmma
Let  $C$ be any of the coordinates
$V_{a1}$ or $H_{ab}$ where $a>1$.
If $V_{11}(\Gamma^-)$ is finite, then $C(\Gamma^-)$ is also finite.
\elmma
%%%%%%%%%%%%%%%%%%%%%%%%%%%%%%%%%%%%%%%%%
\prf
Assume $K$ is a positive integer such that $V_{11}(\Gamma^-)\seq [0,K]$.
Now suppose, if possible, that $C(\Gamma^-)$ is infinite.
Let $\bldr_n$ be a sequence of points in $\Gamma^-$ such that
\[
C(\bldr_1)<C(\bldr_2)<\ldots.
\]
Start at $\bldr_n$, and use lemma~\ref{sweep2} to reach a point
$\bldr'_n$ and then apply $M_{11}$ for $K+1$ times to get to a point
$\blds_n$ for which $V_{11}(\blds_n)>K$.
Throughout this path, $C(\cdot)$ is constant, so that the paths are all
disjoint.
Since $V_{11}(\blds_n)>K$, we have $\blds_n\in \Gamma^+$.
Thus this gives us an infinite ladder for $(\Gamma^+,\Gamma^-)$,
which is impossible.\qed

%%%%%%%%%%%%%%%%%%%%%%%%%%%%%%%%%%%%%%%%%%%%%%%%%%%%
\blmma
Suppose $H_{1\ell}(F)$ is bounded.
If $V_{11}(\Gamma^-)$ is finite, then $F^-$ is  finite.
\elmma
%%%%%%%%%%%%%%%%%%%%%%%%%%%%%%%%%%%%%%%%%%%%%%%%%%%%
\prf
The previous lemma, along with the assumption here
tells us that the sets $V_{a1}(F^-$) and $H_{a,\ell+1-a}(F^-)$
are all bounded for $1\leq a\leq\ell$.
Since for an $\bldr\in V$, one has
$r_{11}=\sum_{a=1}^\ell V_{a1}(\bldr)+\sum_{a=1}^\ell H_{a,\ell+1-a}(\bldr)$,
the set
$\{r_{11}:\bldr\in F^-\}$ is bounded.
It follows that $F^-$ is finite.\qed
%%%%%%%%%%%%%%%%%%%%%%%%%%%%%%%%%%%%%%%%%%%%%%%%%%%%
\bcrlre\label{signcomp}
If $V_{11}(\Gamma^-)$ is finite, then $\scrc ^-$ is  finite.
\ecrlre
%%%%%%%%%%%%%%%%%%%%%%%%%%%%%%%%%%%%%%%%%%%%%%%%%%%%
\prf
Follows from the observation that $H_{1\ell}(\bldr)=0$
for all $\bldr\in \scrc $.\qed\\
A similar argument will tell us that
if $V_{11}(\Gamma^+)$ is finite, then $\scrc ^+$ is  finite.
Thus from lemma~\ref{about_v11}, it follows that
either $\scrc^+$ or $\scrc^-$ is finite.

%%%%%%%%%%%%%%%%%%%%%%%%%%%%%%%%%%%%%%%%%%%%%%%%%%%%%%%%%%%
\bthm\label{eqsign}
Let $D$ be an equivariant Dirac operator on $L_2(SU_q(\ell+1))$.
Then $\sgn D$ must be of the form $2P-I$
or $I-2P$ where $P$ is, up to a compact perturbation, the projection
onto the closed span of
$\{e^\lambda_{\bldr,\blds}: \bldr\in \scrf_{\bldr_i} \mbox{ for some }i\}$,
with $\bldr_1,\ldots,\bldr_k$ being a finite collection of GT-tableaux.
\ethm
%%%%%%%%%%%%%%%%%%%%%%%%%%%%%%%%%%%%%%%%%%%%%%%%%%%%%%%%%%%
\prf
Let
$\scrc '=\{\bldr\in \scrc : \scrf_\bldr^+ \neq\phi\neq \scrf_\bldr^-\}$.
Let us first show that $\scrc '$ is finite, i.e.\ except for
finitely many $\bldr$'s in $\scrc $, one has either
$\scrf_\bldr\seq \Gamma^+$ or $\scrf_\bldr\seq \Gamma^-$.
It follows from the argument used in the proof of theorem~\ref{compl}
that any two points on a free plane can be connected by a
path lying entirely on the plane. If $\scrc'$ is infinite,
one can easily produce an infinite ladder using this fact.

Thus there are only finitely many free
planes $\scrf_\bldr$ for which
both $\scrf_\bldr^+$ and $\scrf_\bldr^-$ are nonempty.
Since we already know that for every $\bldr$,
either  $\scrf_\bldr^+$ or $\scrf_\bldr^-$ is finite,
it follows that by applying a compact perturbation,
one can ensure that for every $\bldr$, exactly one of the sets
$\scrf_\bldr^+$ and $\scrf_\bldr^-$ is empty.
This, along with the observations that
$\scrc\cap\scrf_\bldr=\{\bldr\}$ and that
either $\scrc^+$ or $\scrc^-$ is finite
gives us the required conclusion.\qed

As a consequence of this sign characterization, we now get
the following theorem.

%%%%%%%%%%%%%%%%%%%%%%%%%%%%%%%%%%%%%%%%%%%%%%%
\bthm
Let $\ell>1$. Let $D$ be an  equivariant
Dirac operator acting on $L_2(G)$. Then
$D$ must have trivial sign.
\ethm
%%%%%%%%%%%%%%%%%%%%%%%%%%%%%%%%%%%%%%%%%%%%%%%
\prf
We will show that if $P$ is as in the earlier theorem,
then the commutators $[P,\pi(u_{ij})]$ can not all be compact.

Let us first prove it
 in the case when $P$ is the projection onto the span of
$\{e_{\bldr\blds}: \bldr\in\scrf_0\}$, where $\scrf_0$
is the free plane passing through
the constant tableaux.
We have
\[
[P,\pi(u_{ij})]e_{\bldr\blds}=\cases{
     P\pi(u_{ij})e_{\bldr\blds} & if $\bldr\not\in \scrf_0$,\cr
     (P-I)\pi(u_{ij})e_{\bldr\blds} & if $\bldr\in \scrf_0$}.
\]
Recall (section~5) the expression for $\pi(u_{ij})e_{\bldr\blds}$:
\[
\pi(u_{ij})e_{\bldr\blds}
=\sum_{{R\in\bbn^i, S\in\bbn^j}\atop{R(1)=S(1)}}
     C_q(i,\bldr,R(\bldr))C_q(j,\blds,S(\blds))k(\bldr,R(\bldr))
 e_{R(\bldr)S(\blds)}.
\]
Hence for $\bldr\in \scrf_0$,
\bean
[P,\pi(u_{ij})]e_{\bldr\blds} &=& (P-I)\pi(u_{ij})e_{\bldr\blds}\\
&=& -\sum_{{R\in\bbn^i, S\in\bbn^j}\atop{R(1)=S(1),R\neq N_{i0}}}
 C_q(i,\bldr,R(\bldr))C_q(j,\blds, S(\blds))k(\bldr,R(\bldr))
     e_{R(\bldr),S(\blds)}.
\eean
In particular, for $i=j=1$, one gets
\[
[P,\pi(u_{11})]e_{\bldr\blds} =
 -\sum_{k=1}^\ell
 C_q(1,\bldr,M_{k1}(\bldr))C_q(1,\blds, M_{k1}(\blds))k(\bldr,M_{k1}(\bldr))
     e_{M_{k1}(\bldr),M_{k1}(\blds)}.
\]
Now suppose $\bldr\in\scrf_0$ satisfies
\be \label{choices}
 r_{1,\ell}=0=r_{2,\ell}=r_{1,\ell+1}.
\ee
Then
\[
\langle e_{M_{\ell 1}(\bldr),M_{\ell 1}(\bldr)},
           [P,\pi(u_{11})]e_{\bldr\bldr}\rangle =
 - C_q(1,\bldr, M_{\ell 1}(\bldr))^2k(\bldr,M_{\ell 1}(\bldr)).
\]
It follows from~(\ref{cgc4}) and (\ref{cgc5})
that $C_q(1,\bldr, M_{\ell 1}(\bldr))$ is bounded away from zero,
so long as $\bldr$ obeys (\ref{choices}).
We have also seen (lemma~\ref{krmbound}) that
$k(\bldr,M_{\ell 1}(\bldr))$ is bounded away from zero.
Now it is easy to see that if $\ell>1$, then there are infinitely many choices
of $\bldr$ satisfying (\ref{choices}) such that they all lie in $\scrf_0$.
Therefore $[P,\pi(u_{11})]$
is not compact.

For more general $P$ (as in the previous theorem),
the idea would be similar, but this time
one has to get hold of a positive integer $n$ such that
for any $\bldr\in\cup_{i=1}^k\scrf_{\bldr_i}$,
$nM_{\ell 1}(\bldr)\not\in \cup_{i=1}^k\scrf_{\bldr_i}$,
and then compute
$\langle e_{nM_{\ell 1}(\bldr),nM_{\ell 1}(\bldr)},
                 (P-I)\pi(u_{11})^n e_{\bldr \bldr}\rangle$.
\qed

As mentioned in the introduction, the above theorem
in particular says that in order to get equivariant
Dirac operators  with nontrivial
sign for for $\ell>1$,
one needs to bring in multiplicities. We will see below
that if one takes the tensor product of
$L_2(G)$ with a suitable space, it is possible to produce
such operators.

\bthm
Let $\widetilde{D}$ be as in theorem~\ref{singular} and let
$N_i$ be the following  operators on $L_2(G)$:
\[
N_i e_{\bldr,\blds}=f_i(\bldr) e_{\bldr,\blds},
\]
where $f_i(\bldr)=\min\{H_{ai}(\bldr):1\leq a\leq \ell+1-i\}$.
Let $\gamma_1,\gamma_2,\ldots,\gamma_{\ell+1}$ be
 $\ell+1$ spin matrices acting on $\bbc^m$. Define
 an operator $D$ on $L_2(G)\otimes\bbc^m$ as follows:
\[
D = \sum_{i=1}^\ell N_i\otimes \gamma_i
       + \widetilde{D}\otimes \gamma_{\ell+1}.
\]
Then $(L_2(G)\otimes\bbc^m,\pi\otimes I, D)$ is
an equivariant
$\ell(\ell+2)$-summable spectral triple.

Moreover, the operator $D$ is optimal,
in the following sense:
given any equivariant Dirac operator
$D'$ on $L_2(G)\otimes\bbc^m$
there are positive reals $a,b$  such that
$|D'|\leq a+b|D|$.
\ethm
\prf
Compact resolvent condition and summability
of $D$ follow  from
the fact that the operator $|D|$ is given by
$|D| e_{\bldr,\blds}=\lambda_\bldr e_{\bldr,\blds}$,
where the singular values $\lambda_\bldr$ obey the inequality
\[
r_{11}\leq \lambda_\bldr \leq K r_{11}
\]
for some constant $K$ that depends only on $\ell$.
Boundedness of commutators follow from the boundedness
of commutators of the $N_i$'s and $\widetilde{D}$ with
the algebra elements, which is clear from
condition~(\ref{eqbdd4}).

Observe that  $\widetilde{D}\otimes I\leq |D|$.
Therefore optimality follows from lemma~\ref{optimality1}.
\qed

\brmrk
Let $\widehat{V}_{i1}$ and $\widehat{H}_{ij}$ denote the following
operators on $L_2(G)$:
\[
\widehat{V}_{i1}e_{\bldr,\blds}=V_{i1}(\bldr)e_{\bldr,\blds},
\quad
\widehat{H}_{ij}e_{\bldr,\blds}=H_{ij}(\bldr)e_{\bldr,\blds},
\quad
i+j\leq \ell+1.
\]
Suppose now that $\gamma_1,\gamma_2,\ldots,\gamma_{\ell(\ell+3)/2}$
be spin matrices acting on some space $\bbc^m$, and
$D_k$ for $1\leq k\leq \frac{\ell(\ell+3)}{2}$ are the operators
$\widehat{V}_{i1}$ and $\widehat{H}_{ij}$ in some order.
Now define $D$ on $L_2(G)\otimes\bbc^m$  to be the operator
\[
D=\sum D_k\otimes \gamma_k.
\]
Then this operator $D$ also enjoys all the features
described in the above theorem.
\ermrk

%%%%%%%%%%%%%%%%%%%%%%%%%%%%%%%%%%%%%%%%%%%%
%%%%%%%%%%%%%%%%%%%%%%%%%%%%%%%%%%%%%%%%%%%%
%%%%%%%%%%%%%%%%%%%%%%%%%%%%%%%%%%%%%%%%%%%%
\section{The odd dimensional quantum spheres}
%%%%%%%%%%%%%%%%%%%%%%%%%%%%%%%%%%%%%%%%%%%%
%%%%%%%%%%%%%%%%%%%%%%%%%%%%%%%%%%%%%%%%%%%%
In this section, we will use the combinatorial
technique and the calculations done in the earlier
sections to investigate equivariant Dirac operators
for all the odd dimensional quantum spheres $S_q^{2\ell+1}$
of Vaksman \& Soibelman~(\cite{v-s}).
In what follows, we will write $G$ for $SU_q(\ell+1)$
and  $H$ for $SU_q(\ell)$.

The $C^*$-algebra $C(S_q^{2\ell+1})$ of the quantum
sphere $S_q^{2\ell+1}$
is the universal $C^*$-algebra generated by
elements
$z_1, z_2,\ldots, z_{\ell+1}$
satisfying the following relations (see~\cite{h-s}):
\bean
z_i z_j & =& qz_j z_i,\qquad 1\leq j<i\leq \ell+1,\\
z_i z_j^* & =& q z_j^* z_i ,\qquad 1\leq i\neq j\leq \ell+1,\\
z_i z_i^* - z_i^* z_i +
(1-q^{2})\sum_{k>i} z_k z_k^* &=& 0,\qquad \hspace{2em}1\leq i\leq \ell+1,\\
\sum_{i=1}^{\ell+1} z_i z_i^* &=& 1.
\eean
Just like their classical counterparts,
these spheres can be viewed as quotient spaces
of the quantum groups $SU_q(\ell+1)$, i.\ e.\
\be
C(S_q^{2\ell+1}) \cong C(G\verb1\1H) =
   \{a\in C(G): (\phi\otimes id)\Delta (a)=I\otimes a\},
\ee
where $\phi$ is a $C^*$-homomorphism
from $C(G)$ onto $C(H)$ that preserves the comultiplication,
that is, it satisfies
$\Delta\phi=(\phi\otimes\phi)\Delta$, where the $\Delta$
on the right hand side is the comultiplication for
$G$ and the $\Delta$ on the left hand side stands for the
 comultiplication for $H$.
(For a formulation of quotient spaces etc.\ in the context
of compact quantum groups, see~\cite{po})

The group $G$ has a canonical right action
$\tau:C(G\verb1\1H)\rightarrow C(G\verb1\1H)\otimes C(G)$
coming from the comultiplication $\Delta$
(i.\ e.\ $\tau$ is just the restriction of $\Delta$ to $C(G\verb1\1H)$).
Let $\rho$ denote the restriction of the Haar state on $C(G)$
to $C(G\verb1\1H)$.
Then clearly one has
$(\rho\otimes id)\tau (a) = \rho(a)I$,
which means $\rho$ is the invariant state for $C(G\verb1\1H)$.
This also means that $L_2(G\verb1\1H)=L_2(\rho)$ is just the
closure of $C(G\verb1\1H)$ in $L_2(G)$.

\bppsn
Assume $\ell>1$.
The right regular representation $u$ of $G$ keeps
$L_2(G\verb1\1H)$ invariant, and the restriction of $u$ to
$L_2(G\verb1\1H)$ decomposes as a direct sum of exactly one copy
of each of the irreducibles given by the young tableaux
$\lambda_{n,k}:=(n+k, k,k,\ldots, k,0)$, with $n,k\in\bbn$.
\eppsn
\prf
Write $\sigma$ for the composition $h_H\circ\phi$
where $h_H$ is the Haar state for $H$.
From the description of $C(G\verb1\1H)$ above, it follows
that
\bean
C(G\verb1\1H) &=& \{a\in C(G): (\sigma \otimes id)\Delta (a)=a\}\\
& =&\{(\sigma \otimes id)\Delta (a): a\in C(G)\}.
\eean
Now
the map $a\mapsto \sigma\ast a:=(\sigma\otimes id)\Delta(a)$
on $C(G)$ extends to a bounded linear operator $L_\sigma$ on $L_2(G)$
(lemma~3.1, \cite{pa}), and it is easy to see that
$L_\sigma^2=L_\sigma$. It follows then that
$L_2(G\verb1\1H)=\ker(L_\sigma -I)=\mbox{ran}\,L_\sigma$.
From the discussion preceeding theorem~3.3, \cite{pa},
it now follows that $u$ keeps $L_2(G\verb1\1H)$ invariant and in fact
the restriction of $u$ to $L_2(G\verb1\1H)$ is the representation
induced by the trivial repersentation of $H$.
From the analogue of Frobenius reciprocity theorem for
compact quantum groups (theorem~3.3, \cite{pa}) it now
follows that the multiplicity of any irreducible $u^\lambda$
in it would be same as the multiplicity of the trivial
representation of $H$ in the restriction of $u^\lambda$ to $H$.
But from the representation theory of $SU_q(\ell+1)$,
we know that  the restriction of $u^\lambda$ to $SU_q(\ell)$
decomposes into a direct sum of one copy of
each irreducible $\mu:(\mu_1\geq \mu_2\geq \ldots \geq\mu_\ell)$
of $SU_q(\ell)$ for which
\be\label{induced}
\lambda_1\geq \mu_1 \geq \lambda_2\geq \mu_2\geq \ldots
         \geq\lambda_\ell \geq \mu_\ell \geq 0.
\ee
Now the trivial representation of $SU_q(\ell)$ is indexed
by Young tableaux of the form
$\mu:(k,k,\ldots,k)$ where $k\in\bbn$.
But such a $\mu$ will obey the restriction~\ref{induced} above
if and only if $\lambda$ is of the form
$(n+k,k,k,\ldots,k,0)$.
\qed
\brmrk
For the case $\ell=1$, the restriction of the irreducible
$(n,0)$ to the trivial subgroup decomposes into $n+1$ copies
of the trivial  representation. Therefore, in this case,
$L_2(S_q^3)$ decomposes into a direct sum of $n+1$ copies of
each representation $(n,0)$.
\ermrk

Next, we will make an explicit choice of $\phi$
that would help us make use of the calculations
already done in the initial sections for analyzing
Dirac operators acting on $L_2(G\verb1\1H)$.
More specifically, we will choose our $\phi$ in such
a manner that $L_2(G\verb1\1H)$ turns out to be
the span of certain rows of the $e_{\bldr,\blds}$'s.
Let $u^\one$ denote the fundamental unitary for $G$,
i.\ e.\ the irreducible unitary representation corresponding to the
Young tableaux $\one=(1,0,\ldots,0)$.
Similarly write $v^\one$ for the fundamental unitary for $H$.
Fix some bases for the corresponding representation spaces.
Then $C(G)$ is the $C^*$-algebra generated by the matrix
entries $\{u^\one_{ij}\}$ and $C(H)$ is the
$C^*$-algebra generated by the matrix
entries $\{v^\one_{ij}\}$.
Now define $\phi$ by
\be
\phi(u^\one_{ij})=\cases{ I & if $i=j=1$,\cr
        v^\one_{i-1,j-1} & if $2\leq i,j\leq \ell+1$,\cr
        0 & otherwise.}
\ee
Then $C(G\verb1\1H)$ is the $C^*$-subalgebra of
$C(G)$ generated by the entries $u_{1,j}$ for $1\leq j\leq \ell+1$
(one recovers the relations for the generators of $C(S_q^{2\ell+1})$
if one sets $z_i=q^{-i+1}u^*_{1,i})$.

\bppsn
Let $\Gamma_0$ be the set of all GT tableaux $\bldr^{nk}$
given by
\[
r^{nk}_{ij}=\cases{ n+k & if $i=j=1$,\cr
                0  & if $i=1$, $j=\ell+1$,\cr
                k  & otherwise,}
\]
for some $n,k \in \bbn$.
Let $\Gamma_0^{nk}$ be the set of all GT tableaux with
top row $(n+k,k,\ldots,k,0)$.
Then the family of vectors
\[
\{e_{\bldr^{nk},\blds}: n,k\in\bbn,\, \blds\in\Gamma_0^{nk}\}
\]
form a complete
orthonormal basis for $L_2(G\verb1\1H)$.
\eppsn
\prf
Let $A$ be the linear span of the elements
$\{u_{\bldr^{n,k},\blds}: n,k\in\bbn, \blds\in\Gamma_0^{n,k}\}$.
Clearly  the closure of $A$ in $L_2(G)$ is the closed
linear span of $\{e_{\bldr^{nk},\blds}: n,k\in\bbn,\, \blds\in\Gamma_0^{nk}\}$.
It is also immdiate that the restriction of the
right regular representation to the above subspace
is a direct sum of one copy of each of the irreducibles
$(n+k,k,k,\ldots,k,0)$.

We will next show that for any $a\in A$, $u_{1j}a$ and $u_{1j}^*$ a are also
in $A$.
Take $a=u_{\bldr^{n,k},\blds}$. Use equation~(\ref{alg_left_mult}) to get
\bea
u_{1,j}u_{\bldr^{n,k},\blds} &=&
  \sum_{M, M'}
  C_q(1,\bldr^{n,k},M(\bldr^{n,k}))C_q(j,\blds,M'(\blds))
      u_{M(\bldr^{n,k}),M'(\blds)}\cr
 &=& \sum_{M'}
  C_q(1,\bldr^{n,k},M_{11}(\bldr^{n,k}))C_q(j,\blds,M'(\blds))
      u_{M_{11}(\bldr^{n,k}),M'(\blds)} \cr
&&   +
   \sum_{M''}
  C_q(1,\bldr^{n,k},M_{\ell+1,1}(\bldr^{n,k}))C_q(j,\blds,M''(\blds))
      u_{M_{\ell+1,1}(\bldr^{n,k}),M''(\blds)} \cr
&=&\sum_{M'}
  C_q(1,\bldr^{n,k},\bldr^{n+1,k})C_q(j,\blds,M'(\blds))
      u_{\bldr^{n+1,k},M'(\blds)}\cr
&&   +
   \sum_{M''}
  C_q(1,\bldr^{n,k},\bldr^{n,k-1}))C_q(j,\blds,M''(\blds))
      u_{\bldr^{n,k-1},M''(\blds)},
\eea
where the first sum is over all moves $M'\in\bbn^{j}$
whose first coordinate is 1 and the second
sum is over all moves $M''\in\bbn^{j}$
whose first coordinate is $\ell+1$.
Thus $u_{1j}a\in A$.

Next, note that if
$\langle u_{1j}^* e_{\bldr^{n,k},\blds}, e_{\bldr',\blds'}\rangle\neq 0$,
then one must have
$\bldr'=\bldr^{n-1,k}$ or $\bldr'=\bldr^{n,k+1}$.
Therefore it follows that $u_{1j}^* u_{\bldr^{n,k},\blds}$
is a linear combination of the $u_{\bldr^{n-1,k},\blds}$
$u_{\bldr^{n,k+1},\blds}$'s, and hence belongs to $A$.
Since $A$ contains the element $u_{\mathbf{0},\mathbf{0}}=1$,
it contains $u_{1j}$ and $u_{ij}^*$. Thus $A$ contains the $*$-algebra
$B$ generated by the $u_{1j}$'s.
But by the previous theorem, restriction of the right regular representation
to the $L_2$ closure $L_2(G\verb1\1H)$ of $B$ also decomposes
as a direct sum of one copy
of each of the irreducibles $(n+k,k,\ldots,k,0)$.
So it follows that $L_2(G\verb1\1H)$ is equal to
the subspace stated in the theorem.
\qed

A self-adjoint operator with compact resolvent on $L_2(G\verb1\1H)$
that commutes with the restriction of $u$ there would be
of the form
\[
e_{\bldr,\blds}\mapsto d(\bldr)e_{\bldr,\blds},\quad \bldr\in\Gamma_0.
\]
Next, let us look at the growth restrictions coming from the
boundedness of commutators.
In this case, one has the boundedness of only the operators
$[D,\pi(u_{ij})]$. Which means, in effect, one will now have
the condition~(\ref{eqbdd4}) only for $i=1$ and $\bldr\in\Gamma_0$:
\be\label{eqbdd_sph1}
|d(\bldr)-d(M(\bldr))|\leq c q^{-C(1,\bldr,M)}.
\ee
Observe that only allowed moves here are the moves
$M=M_{1,1}\equiv(1)$ and $M=M_{\ell+1,1}\equiv(\ell+1)$.
Looking at the corresponding quantity
$C(1,\bldr,M)$, we find that there are two conditions:
\bea
|d(\bldr^{nk})-d(\bldr^{n,k-1})| &\leq & c,\label{eqbdd_sph2}\\
|d(\bldr^{nk})-d(\bldr^{n+1,k})| &\leq &
       cq^{-\sum_{j=1}^{\ell}H_{1j}(\bldr^{nk})}
      =cq^{-k}.\label{eqbdd_sph3}
\eea
As in the earlier sections, we can now
form a graph by taking $\Gamma_0$ to be the set of vertices,
and by joining two vertices $\bldr$ and $\blds$ by an edge if
$|d(\bldr)-d(\blds)|\leq c$.

\blmma
Let $\scrf_n=\{\bldr^{n,k}:k\in\bbn\}$, $n\in\bbn$.
Then any two points in $\scrf_n$ are connected
by a path lying entirely in $\scrf_n$.

If $n<n'$, then any point in $\scrf_n$ is connected to
any point in $\scrf_{n'}$ by a path such that
$n\leq V_{1,1}(\bldr) \leq n'$
for every vertex $\bldr$ lying on that path.
\elmma
\prf
Take two points $\bldr^{n,j}$ and $\bldr^{n,k}$
in $\scrf_n$. Assume $j<k$.
From the condition (\ref{eqbdd_sph2}), it follows
that any point $\bldr$ is connected to $M_{\ell+1,1}(\bldr)$
by an edge. Therefore the first conclusion follows
from the observation that if we start at $\bldr^{n,k}$
and apply the move $M_{\ell+1,1}$ successively $k-j$ number of times,
we reach the point $\bldr^{n,j}$, and the vertices on this path
are the points $\bldr^{n,i}$ for $i=j, j+1,\ldots,k$.
Observe also that throughout this path, $V_{1,1}(\bldr)$
remains $n$.

For the second part, take a point $\bldr^{n,k}$ in $\scrf_n$
and a point $\bldr^{n',j}$ in $\scrf_{n'}$.
From what we have done above, there is a path
from $\bldr^{n,k}$ to $\bldr^{n,0}$ throughout which
$V_{1,1}(\bldr)=n$.
Similarly there is a path
from $\bldr^{n',j}$ to $\bldr^{n',0}$ throughout which
$V_{1,1}(\bldr)=n'$.
Next, note from (\ref{eqbdd_sph3}) that  for $p\in\bbn$,
the points
$\bldr^{p,0}$ and $\bldr^{p+1,0}$ are connected by an edge
and
$V_{1,1}(\bldr^{p,0})=p$, $V_{1,1}(\bldr^{p+1,0})=p+1$.
So start at $\bldr^{n,0}$ and reach successively the
points
$\bldr^{n+1,0}$, $\bldr^{n+2,0}$ and so on to
eventually reach the point $\bldr^{n',0}$;
also the coordinate $V_{1,1}(\cdot)$ remains between $n$ and $n'$
on this path.\qed

\bthm\label{eqsign_sphere}
Let $D$ be an equivariant Dirac operator on $L_2(G\verb1\1H)$.
Then
\begin{enumerate}
\item
$D$ must be of the form
\[
e_{\bldr,\blds}\mapsto d(\bldr)e_{\bldr,\blds},\quad \bldr\in\Gamma,
\]
where the singular values obey $|d(\bldr)|=O(r_{11})$, and
\item
$\sgn D$ must be of the form $2P-I$
or $I-2P$ where $P$ is, up to a compact perturbation, the projection
onto the closed span of
$\{e_{\bldr^{nk},\blds}: n\in F, k\in\bbn, \blds\in \Gamma_0^{nk}\}$,
for some finite subset $F$ of $\bbn$.
\end{enumerate}
\ethm
\prf
Start with an equivariant self-adjoint operator
$D$ with compact resolvent, so that it is indeed of the form
$e_{\bldr,\blds}\mapsto d(\bldr)e_{\bldr,\blds}$.
By applying a compact perturbation if necessary,
make sure that $d(\bldr)\neq 0$ for all $\bldr\in\Gamma_0$.
We have seen during the proof of the previous lemma that
for any $n$ and $k$ in $\bbn$, the vertices
$\bldr^{nk}$ and $\bldr^{n,k+1}$ are connected by an edge,
and for any $n\in\bbn$, the vertices
$\bldr^{n,0}$ and $\bldr^{n+1,0}$ is connected by an edge.
Thus any vertex $\bldr^{nk}$ can be reached from the vertex
$\bldr^{00}$ by a path of length $n+k$. Therefore one gets the
first assertion.

Next, define
\bean
\Gamma_0^+ &=& \{\bldr\in\Gamma_0: d(\bldr)>0\},\\
\Gamma_0^- &=& \{\bldr\in\Gamma_0: d(\bldr)<0\},\\
\scrf_n^+ &=& \scrf_n\cap \Gamma_0^+,\\
\scrf_n^- &=&  \scrf_n\cap \Gamma_0^-.
\eean
Observe that for the path produced in the proof
of the forgoing lemma to connect two
points $\bldr^{n,k}$ and $\bldr^{n,j}$ in $\scrf_n$,
the coordinate $H_{1,\ell}(\cdot)$ remains between $j$ and $k$.
Now suppose for some  $n$,
both $\scrf_n^+$ and $\scrf_n^-$ are infinite.
Then there are points
\[
0\leq k_1 < k_2 < \ldots
\]
such that $\bldr^{nk}$ is in $\scrf_n^+$ for $k=k_{2j}$
and  $\bldr^{nk}$ is in $\scrf_n^-$ for $k=k_{2j+1}$.
Using the above observation, we can then produce
an infinite ladder by joining
each $\bldr^{n,k_{2j-1}}$ to $\bldr^{n,k_{2j}}$.
Thus for each $n\in\bbn$, exactly one of the sets
$\scrf_n^+$ and $\scrf_n^-$ is finite.
Also, note that by the first part of the previous lemma,
the set of all $n\in\bbn$ for which
both $\scrf_n^+$ and $\scrf_n^-$ are nonempty is finite.
Therefore by applying a compact perturbation, we can ensure that
for every $n$, either $\scrf_n^+=\scrf_n$ or
$\scrf_n^-=\scrf_n$.

Finally, if there are infinitely many $n$'s for which
$\scrf_n^+=\scrf_n$
and infinitely many $n$'s for which $\scrf_n^-=\scrf_n$,
then one can choose a sequence of integers
\[
0\leq n_1 < n_2 <\ldots
\]
such that
$\scrf_n^+=\scrf_n$ for $n=n_{2j}$
and
$\scrf_n^-=\scrf_n$ for $n=n_{2j+1}$.
Now use the second part of the previous lemma
to join each $\bldr^{n_{2j-1},0}$ to $\bldr^{n_{2j},0}$
to produce an infinite ladder.

Thus there is a finite subset $F$ of $\bbn$ such that
exactly one of the following is true:
\[
\scrf_n=\cases{\scrf_n^+ & if $n\in F$,\cr
               \scrf_n^- & if $n\not\in F$,}
\qquad
\mbox{or }
\qquad
\scrf_n=\cases{\scrf_n^- & if $n\in F$,\cr
               \scrf_n^+ & if $n\not\in F$.}
\]
This is precisely what the second part
of the theorem says.\qed

Next, take the operator $D:e_{\bldr,\blds}\mapsto d(\bldr)e_{\bldr,\blds}$
 on $L_2(G\verb1\1H)$ where the $d(\bldr)$'s are given by:
\be\label{eq_sphere1}
d(\bldr^{nk})=\cases{-k & if $n=0$,\cr
                     n+k & if $n>0$.}
\ee
\bthm\label{generic_d_sph}
The operator $D$ is an equivariant $(2\ell+1)$-summable
Dirac operator acting on $L_2(G\verb1\1H)$, that gives a
nondegenerate pairing with $K_1(C(G\verb1\1H))$.

The operator $D$ is optimal, i.\ e.\
if $D_0$ is any equivariant Dirac operator on $L_2(G\verb1\1H)$,
then there are positive reals $a$ and $b$ such that
\[
|D_0|\leq a+b|D|.
\]
\ethm
\prf
Recall from equation~(\ref{left_mult}) that the elements
$u_{1,j}$ act on the basis
elements $e_{\bldr^{n,k},\blds}$ as follows:
\bea\label{l2_repn_sph}
u_{1,j}e_{\bldr^{n,k},\blds} &=&
  \sum_{M, M'}
  C_q(1,\bldr^{n,k},M(\bldr^{n,k}))C_q(j,\blds,M'(\blds))
   \kappa(\bldr^{n,k},\blds)
   e_{M(\bldr^{n,k}),M'(\blds)}\cr
 &=& \sum_{M'}
  C_q(1,\bldr^{n,k},M_{11}(\bldr^{n,k}))C_q(j,\blds,M'(\blds))
   \kappa(\bldr^{n,k},\blds)
   e_{M_{11}(\bldr^{n,k}),M'(\blds)} \cr
&&   +
   \sum_{M''}
  C_q(1,\bldr^{n,k},M_{\ell+1,1}(\bldr^{n,k}))C_q(j,\blds,M''(\blds))
   \kappa(\bldr^{n,k},\blds)
   e_{M_{\ell+1,1}(\bldr^{n,k}),M''(\blds)} \cr
&=&\sum_{M'}
  C_q(1,\bldr^{n,k},\bldr^{n+1,k})C_q(j,\blds,M'(\blds))
   \kappa(\bldr^{n,k},\blds)
   e_{\bldr^{n+1,k},M'(\blds)}\cr
&&   +
   \sum_{M''}
  C_q(1,\bldr^{n,k},\bldr^{n,k-1}))C_q(j,\blds,M''(\blds))
   \kappa(\bldr^{n,k},\blds)
   e_{\bldr^{n,k-1},M''(\blds)},
\eea
where the first sum is over all moves $M'\in\bbn^{j}$
whose first coordinate is 1 and the second
sum is over all moves $M''\in\bbn^{j}$
whose first coordinate is $\ell+1$.
If we now plug in the values of the Clebsch-Gordon coefficients
from equations~(\ref{cgc4}) and~(\ref{cgc5}), we get
\bea
u_{1,j}e_{\bldr^{n,k},\blds} &=&
 \sum_{M'}
  P'_1 P'_2 q^{k+C(j,\blds,M')}
   \kappa(\bldr^{n,k},\blds)
   e_{\bldr^{n+1,k},M'(\blds)}\cr
&&     +
   \sum_{M''}
  P''_1 P''_2 q^{C(j,\blds,M'')}
   \kappa(\bldr^{n,k},\blds)
   e_{\bldr^{n,k-1},M''(\blds)},
\eea
where $P'_i$, $P''_j$ and $k(\bldr^{n,k},\blds)$
all lie between two fixed positive numbers.
Boundedness of the commutators $[D,u_{1,j}]$
now follow directly.

For summability, notice that the eigenspace
of $|D|$ corresponding to the eigenvalue $n\in\bbn$
is the span of
\[
\{e_{\bldr^{k,n-k},\blds}: 0\leq k\leq n, \blds\in\Gamma_0^{k,n-k}\}.
\]
Now just count the number of elements in the above set
to get summability.

Next, we will compute the pairing of the $K$-homology class
of this $D$ with a generator of the $K_1$ group.
Write $\omega_q:=q^{-\ell}u_{1,\ell+1}$.
From the commutation relations, it follows that
this element has spectrum
\[
\{z\in\bbc: |z|=0 \mbox{ or }q^n\mbox{ for some }n\in\bbn\}.
\]
Then the element $\gamma_q:=\chi_{\{1\}}(\omega_q^*\omega_q)(\omega_q-I)+I$
is unitary.
We will show that the index of the operator
$Q\gamma_q Q$ (viewed as an operator on $QL_2(G\verb1\1H)$)
is $1$, where $Q=\frac{I-\sgn D}{2}$, i.\ e.\ it
is the projection onto the closed
linear span of
$\{e_{\bldr^{0,k},\blds}:k\in\bbn, \blds\in\Gamma_0^{0,k}\}$.
What we will actually do is
compute the index of the
operator $Q\gamma_0 Q$ and appeal to continuity of the index.
From equation~(\ref{l2_repn_sph}), we get
\bea\label{for_q=0_1}
\lefteqn{u_{1,\ell+1}e_{\bldr^{0,k},\blds}}\cr
 &=&
  C_q(1,\bldr^{0,k},M_{11}(\bldr^{0,k}))C_q(\ell+1,\blds,N_{1,0}(\blds))
    \kappa(\bldr^{0,k},M_{11}(\bldr^{0,k}))e_{\bldr^{1,k},N_{1,0}(\blds)}\cr
&& +     C_q(1,\bldr^{0,k},M_{\ell+1,1}(\bldr^{0,k}))
                C_q(\ell+1,\blds,M_{\ell+1,\ell+1}(\blds))
    \kappa(\bldr^{0,k},M_{\ell+1,1}(\bldr^{0,k}))
        e_{\bldr^{0,k-1},M_{\ell+1,\ell+1}(\blds)}.\cr
&&
\eea
Use the formula~(\ref{cgc1}) for Clebsch-Gordon coefficients
to get
\bea
C_q(1,\bldr^{0,k},M_{11}(\bldr^{0,k}))
 &=& q^{k}(1+ o(q)),\\
C_q(1,\bldr^{0,k},M_{\ell+1,1}(\bldr^{0,k})
 &=&  1+ o(q),\\
C_q(\ell+1,\blds,N_{1,0}(\blds))
 &=& 1+ o(q),\\
C_q(\ell+1,\blds,M_{\ell+1,\ell+1}(\blds))
 &=& q^{s_{\ell+1,1}+\ell}(1+ o_4(q)),
\eea
where $o(q)$ signifies a function of $q$ that is
continuous at $q=0$ and
$o(0)=0$.
We also have
\bea
\kappa(\bldr^{0,k},M_{11}(\bldr^{0,k}))
  &=& q^\ell (1+o(q)),\\
\kappa(\bldr^{0,k},M_{\ell+1,1}(\bldr^{0,k}))
  &=& 1+o(q),
\eea
where $o(q)$ is as earlier.
Plugging these values in~(\ref{for_q=0_1}) %% and putting $q=0$,
we get
% \be
% \omega_0 e_{\bldr^{0,k},\blds}
%   = \cases{
%    e_{\bldr^{0,k-1},M_{\ell+1,\ell+1}(\blds)}
%              & if $k>0$ and $s_{\ell+1,1}=0$,\cr
%     e_{\bldr^{1,0},N_{1,0}(\blds)} & if $k=0$,\cr
%     0 & otherwise.}
% \ee
\be
\omega_q e_{\bldr^{0,k},\blds}
  =  q^{k}(1+o(q))e_{\bldr^{1,k},N_{1,0}(\blds)}
   + q^{s_{\ell+1,1}}(1+o(q))e_{\bldr^{0,k-1},M_{\ell+1,\ell+1}(\blds)}
\ee
Putting $q=0$,
we get
\be
\omega_0 e_{\bldr^{0,k},\blds}
  = \cases{
   e_{\bldr^{0,k-1},M_{\ell+1,\ell+1}(\blds)}
             & if $k>0$ and $s_{\ell+1,1}=0$,\cr
    e_{\bldr^{1,0},N_{1,0}(\blds)} & if $k=0$,\cr
    0 & otherwise.}
\ee
Thus $\omega_0^*\omega_0$ is the projection onto the span
of
$\{e_{\bldr^{0,k},\blds^k}: k\in\bbn\}$
where
$\blds^k$ is the GT tableaux given by
\[
s^k_{ij}=\cases{0 & if $i=\ell+2-j$,\cr
                  k & otherwise,}
\]
which is uniquely determined by the conditions
$s_{\ell+1,1}=0$ and that $\blds\in\Gamma_0^{0,k}$.
Therefore the operator $\gamma_0$  is given by
\[
\gamma_0 e_{\bldr^{0,k},\blds} =
 e_{\bldr^{0,k},\blds} - \chi_{\{\blds=\blds^k\}}e_{\bldr^{0,k},\blds}
   + \chi_{\{\blds=\blds^k\}}
                   e_{\bldr^{0,k-1},\blds^{k-1}}.
\]
It now follows that the index of $Q\gamma_0 Q$ is $1$.

Optimality follows from part~1 of the previous theorem.
\qed

%%%%%%%%%%%%%%%%%%%%%%%%%%
%%%  BIBLIOGRAPHY
%%%%%%%%%%%%%%%%%%%%%%%%%%

%%%%%%%%%%%%%%%%%%%%%%%%%%%%%%%
\noindent{\sc Partha Sarathi Chakraborty}
(\texttt{chakrabortyps@cf.ac.uk})\\
         {\footnotesize School of Mathematics,
 Cardiff University, Senghennydd Road, Cardiff, UK}\\[1ex]
{\sc Arupkumar Pal} (\texttt{arup@isid.ac.in})\\
         {\footnotesize Indian Statistical
Institute, 7, SJSS Marg, New Delhi--110\,016, INDIA}

%%%%%%%%%%%%%%%%%%%%%%%%%%%%%%%%%%%%%%%%%%%%%%%%%%%%%%%%%%%
%%%%%%%%%%%%%%%%%%%%%%%%%%%%%%%%%%%%%%%%%%%%%%%%%%%%%%%%%%%
%%%%%%%%%%%%%%%%%%%%%%%%%%%%%%%%%%%%%%%%%%%%%%%%%%%%%%%%%%%
%%%%%%%%%%%%%%%%%%%%%%%%%%%%%%%%%%%%%%%%%%%%%%%%%%%%%%%%%%%

\end{document}